\documentclass[12pt, reqno, twoside]{amsart}
\usepackage{amsmath, amsthm, mathrsfs, amscd, amsfonts, amssymb, graphicx, color}
\usepackage[bookmarksnumbered, colorlinks, plainpages]{hyperref}
\hypersetup{colorlinks=true,linkcolor=red, anchorcolor=green, citecolor=cyan, urlcolor=red, filecolor=magenta, pdftoolbar=true}
\usepackage{enumitem}
\usepackage{enumitem}
\newtheorem{theorem}{Theorem}[section]
\newtheorem{lemma}[theorem]{Lemma}

\theoremstyle{definition}

\theoremstyle{remark}

\numberwithin{equation}{section}

\begin{document}

\setcounter{page}{1}

\title[Gradient estimates]{Gradient estimates for a nonlinear parabolic equation on smooth metric measure spaces with evolving metrics and potentials}

\author[A. Taheri]{Ali Taheri}

\author[V. Vahidifar]{Vahideh Vahidifar}

\address{School of Mathematical and  Physical Sciences, 
University of Sussex, Falmer, Brighton, United Kingdom.}
\email{\textcolor[rgb]{0.00,0.00,0.84}{a.taheri@sussex.ac.uk}} 
 
\address{School of Mathematical and  Physical Sciences, 
University of Sussex, Falmer, Brighton, United Kingdom.}
\email{\textcolor[rgb]{0.00,0.00,0.84}{v.vahidifar@sussex.ac.uk}}

\subjclass[2020]{53C44, 58J60, 58J35, 60J60}

\keywords{Smooth metric measure spaces, Witten Laplacian, super Perelman-Ricci flow, Li-Yau estimates, Liouville-type results, Harnack inequalities}


\begin{abstract}
This article presents new parabolic and elliptic type gradient estimates for positive smooth solutions to a nonlinear parabolic 
equation involving the Witten Laplacian in the context of smooth metric measure spaces. The metric and potential here 
are time dependent and evolve under a super Perelman-Ricci flow. The estimates are derived under natural lower 
bounds on the associated generalised Bakry-\'Emery Ricci curvature tensors and are utilised in establishing fairly general 
local and global bounds, Harnack-type inequalities and Liouville-type global constancy theorems to mention a few. Other 
implications and consequences of the results are also discussed.
\end{abstract}

\maketitle
 
{ 
\hypersetup{linkcolor=black}
\tableofcontents
}

\section{Introduction} \label{sec1}

In this paper we study gradient estimates for positive smooth solutions to a class of nonlinear parabolic equations 
on smooth metric measure spaces $(M,g, d\mu)$ with metrics and potentials evolving under a $(\mathsf{k},m)$-super 
Perelman-Ricci flow. More specifically we consider positive smooth solutions 
$u=u(x,t)$ to the coupled system  (for $x \in M$, $t>0$)
\begin{align} \label{eq11}
\begin{cases}
\square_q u(x,t) := \left[ \dfrac{\partial}{\partial t} - q(x,t) - \Delta_f \right] u (x,t) = \Sigma(t,x, u(x,t)), \\
\dfrac{\partial g}{\partial t} (x,t) + 2{\mathscr Ric}^m_f(g)(x,t) \ge - 2\mathsf{k} g(x,t). 
\end{cases}
\end{align}

Here $M$ is a complete Riemannian manifold of dimension 
$n \ge 2$, $d\mu=e^{-f} dv_g$ is a weighted measure associated with the smooth potential $f=f(x,t)$ and $dv_g$ 
is the usual Riemannain volume measure. Note that both the metric $g$ and the potential $f$ are time dependent 
and evolve under a $(\mathsf{k}, m)$-super Perelman-Ricci flow as is formulated by the symmetric $(0,2)$-tensor flow inequality on the 
second line in the system \eqref{eq11}. Referring to the equation on the first line in \eqref{eq11} the operator 
$\Delta_f$ is the Witten Laplacian (also known as the weighted or drifting Laplacian, or at times to emphasise 
the choice of $f$, the $f$-Laplacian) acting on functions $v \in \mathscr{C}^2(M)$ through the identity  
\begin{equation} \label{Lf definition}
\Delta_f v = e^f {\rm div} (e^{-f} \nabla v) = \Delta v - \langle \nabla f, \nabla v \rangle,  
\end{equation}
where $\Delta, {\rm div}$ and $\nabla$ are the usual Laplace-Beltrami, divergence and gradient operators 
associated with the metric $g$ respectively. Moreover, the evolution operators
\begin{equation}
\square_q = \partial_t - q - \Delta_f, \qquad \square = \square_0 = \partial_t - \Delta_f, 
\end{equation} 
are the $q$-weighted (and weighted) heat operators respectively with $q=q(x,t)$ a smooth function of the 
space-time variables.

As indicated above the evolution of the metric-potential pair $(g,f)$ is governed by the flow inequality on the second line in \eqref{eq11} 
with $\mathsf{k}$ in $\mathbb{R}$ and $m \ge n$ (note in particular that $m$ is not necessarily an integer). Now referring to the left-hand 
side of this inequality, the symmetric second order space-time dependent tensor field   
\begin{equation} \label{Ricci-m-f-intro}
{\mathscr Ric}^m_f(g)(x,t) = {\mathscr Ric}(g)(x,t) + {\rm Hess}(f)(x,t) - \frac{\nabla f \otimes \nabla f}{m-n}(x,t), 
\end{equation}
is the generalised Bakry-\'Emery Ricci curvature tensor of the triple $(M,g,d\mu)$ with ${\mathscr Ric}(g)$ the usual Riemannain Ricci 
curvature tensor of $g$ and ${\rm Hess}(f)$ the Hessian of $f$ ({\it see} \cite{BD, BE, Bak}). For the sake of clarity we point out that when 
$m=n$, by convention, $f$ is only allowed to be a constant, resulting in ${\mathscr Ric}^m_f(g)={\mathscr Ric}(g)$, 
whilst, we also allow for $m = \infty$, in which case by formally passing to the limit in \eqref{Ricci-m-f-intro} we set 
\begin{equation} \label{Ricf def eq} 
{\mathscr Ric}_f^\infty(g)(x,t) = {\mathscr Ric}(g)(x,t) + {\rm Hess}(f)(x,t) := {\mathscr Ric}_f(g)(x,t).
\end{equation}

The nonlinearity $\Sigma=\Sigma(t,x,u)$ in \eqref{eq11} is a sufficiently smooth function depending on both the space-time 
variables $(x,t)$ with $x \in M$, $t \ge 0$ and the independent variable $u$. We give several examples of such 
nonlinearities arising from different contexts, e.g., conformal geometry and mathematical physics, each presenting a different 
phenomenon whilst depicting a corresponding singular or regular behaviour in certain variable ranges and regimes.

The $f$-Laplacian \eqref{Lf definition} is a symmetric diffusion operator with respect to the invariant weighted measure 
$d\mu=e^{-f} dv_g$. It arises in a variety of contexts including probability theory and stochastic analysis, geometry, 
quantum field theory, statistical mechanics and kinetic theory \cite{BD, Bak, Gr, VC}. It is a natural generalisation of the Laplace-Beltrami 
operator to the smooth metric measure space context and it coincides with the latter precisely when the potential $f$ is a 
constant. By an application of the integration by parts formula it can been that for $u, w \in {\mathscr C}_0^\infty(M)$ it holds  
\begin{equation}
\int_M e^{-f} w \Delta_f u \, dv_g 
= - \int_M e^{-f} \langle \nabla u, \nabla w \rangle \, dv_g 
= \int_M e^{-f} u \Delta_f w \, dv_g
\end{equation}

Our main objective in this paper is to establish new and fairly general local and global elliptic as well as parabolic type gradient 
estimates for positive smooth solutions to \eqref{eq11} and present some of their implications in a context where the metric 
and potential are time dependent and evolve under a $(\mathsf{k},m)$-super Perelman-Ricci flow.

Gradient estimates occupy a central place in geometric analysis with a huge scope of applications (see  
\cite{[Ba1], Bri, ChYau, Dung, [Ha93],[LY86], [LLi1], [LLi2]}, 
\cite{AM, Bid, CDZV, DungHT, NiL, SZ, Taheri-GE-1, Taheri-GE-2, TVahGrad, TVahGradSZ, [Wang], [Wu15], [Wu18]} 
as well as \cite{Bak, Gr, Giaq, [LiP],Taheri-book-two, [QZh]} and  the references therein). The estimates of interest in this paper fall under 
the category of Souplet-Zhang, Hamilton and Li-Yau types that were first formulated and proved for the linear heat 
equation on static manifolds in \cite{[Ha93], SZ} and later extended by many authors to contexts including evolving manifolds 
and nonlinear equations. To the best of our knowledge the results presented 
here are the first under such generalities on the nonlinearity in \eqref{eq11} and/or subject to the evolution of the 
metric-potential pair. In particular they encompass, unify and extend many existing results in the literature 
for very specific types of nonlinearities.

According to the weighted Bochner-Weitzenb\"ock formula, the Witten Laplacian $\Delta_f$ and the generalised 
curvature tensor ${\mathscr Ric}_f(g)$ are related, for every function $u \in {\mathscr C}^3(M)$, via the identity 
\begin{equation} \label{Bochner}
\frac{1}{2} \Delta_f |\nabla u|^2 = |{\rm Hess}(u)|^2 + \langle \nabla u, \nabla \Delta_f u \rangle + {\mathscr Ric}_f (\nabla u, \nabla u).
\end{equation}

Since by an application of the Cauchy-Schwartz inequality we have $\Delta u \le \sqrt n |{\rm Hess}(u)|$ 
(with the norm of the Hessian on the right being the usual Hilbert-Schmidt norm on $2$-tensors), 
upon recalling $\Delta_f u = \Delta u - \langle \nabla f, \nabla u \rangle$ in \eqref{Lf definition} and 
by using Young's inequality, it is seen that 
\begin {equation} \label{Hessian-Witten}
|{\rm Hess}(u)|^2 + \frac{\langle \nabla f , \nabla u \rangle^2}{m-n} 
\ge \frac{(\Delta u)^2}{n} + \frac{\langle \nabla f, \nabla u \rangle ^2}{m-n} 
\ge \frac{(\Delta u - \langle \nabla f, \nabla u \rangle)^2}{m} = \frac{(\Delta_f u)^2}{m}.  
\end{equation}
Hence \eqref{Hessian-Witten} in conjunction with \eqref{Bochner} gives the inequality 
\begin{equation} \label{Bochner-m-inequality}
\frac{1}{2} \Delta_f |\nabla u|^2 - \langle \nabla u, \nabla \Delta_f u \rangle 
\ge \frac{1}{m} (\Delta_f u)^2  + {\mathscr Ric}^m_f (\nabla u, \nabla u).
\end{equation}
As a result, subject to a curvature lower bound of the form ${\mathscr Ric}_f^m(g) \ge -\mathsf{k} g$ the operator $L=\Delta_f$ 
is seen to satisfy the curvature dimension condition ${\rm CD}(-\mathsf{k},m)$ for symmetric diffusion operators, specifically, 
({\it cf.} \cite{BE, Bak}) 
\begin{equation} \label{Bochner-m-inequality-alt}
\frac{1}{2} L |\nabla u|^2 - \langle \nabla u, \nabla L u \rangle 
\ge \frac{1}{m} (L u)^2  - \mathsf{k} |\nabla u|^2.
\end{equation}
The latter condition and inequality plays a fundamental role in the analysis of diffusion operators and their geometric properties.

Let us recall that in the static case, i.e., non-evolving metric-potential pairs, gradient estimates for positive solutions to linear 
and nonlinear heat type equations have been studied extensively starting from the seminal paper of Li and Yau \cite{[LY86]} 
({\it see} also \cite{[LiP]}). In the nonlinear setting the first equation to be considered is the one with a logarithmic type nonlinearity 
({\it see}, e.g., \cite{Ma})
\begin{equation} \label{eq1.4}
\square_q u = [\partial_t - q(x,t) - \Delta_f] u = p(x,t) u \log u. 
\end{equation}
The interest in such problems originates partly from its natural links with gradient Ricci solitons and partly from links with 
geometric and functional inequalities on manifolds, notably, the logarithmic Sobolev and energy-entropy 
inequalities \cite{BE, Bak, Gross, VC}. Recall that a Riemannian manifold $(M,g)$ is said to be a gradient 
Ricci soliton {\it iff} there there exists a smooth function $f$ on $M$ and a constant $\lambda \in \mathbb{R}$ 
such that ({\it cf.} \cite{Cao, Chow, Lott})
\begin{equation}\label{eq1.5}
{\mathscr Ric}_f(g) = {\mathscr Ric}(g) + {\rm Hess}(f) = \lambda g.
\end{equation}
A gradient Ricci soliton can be expanding ($\lambda>0$), steady ($\lambda=0$) or shrinking ($\lambda<0$). 
The notion is a generalisation of an Einstein manifold and has a fundamental role in the analysis of singularities 
of the Ricci flow \cite{Ham, [QZh]}. Other classes of equations closely relating to \eqref{eq1.4} include 
\begin{equation}
\square_q u = [\partial_t - q(x,t) - \Delta_f] u = a(x,t) \Gamma(\log u) u^\mathsf{p} + b(x,t) u^\mathsf{q}, 
\end{equation}
where $\mathsf{p}, \mathsf{q}$ are real exponents, $a, b$ are sufficiently smooth functions and 
$\Gamma \in {\mathscr C}^2(\mathbb{R}, \mathbb{R})$ ({\it see} \cite{CGS, GKS, Taheri-GE-1, Taheri-GE-2, Wu10, [Wu15]} 
and the references therein for gradient estimates and related results in this direction).

Another class of equations that have been extensively studied and whose nonlinearity is in the form of superposition 
of power-like nonlinearities are the Yamabe equations ({\it see}, e.g., \cite{Bid, Case, GidSp, LiJ91}). In the context of 
smooth metric measure spaces a far-reaching generalisation of these equations take the form ({\it see}, e.g., \cite{[Wu18]}) 
\begin{align} \label{eq11s}
\square_q u = [\partial_t - q(x,t) - \Delta_f] u = a(x,t) u^\mathsf{p} + b(x,t) u^\mathsf{q}. 
\end{align}

A closely related and yet more general form of Yamabe type equations is the Einstein-scalar field Lichnerowicz 
equation (see, e.g., Choquet-Bruhat \cite{CBY}, Chow \cite{Chow} and Zhang \cite{[QZh]}). When the underlying 
manifold has dimension $n \ge 3$ the evolutionary form of this equation takes the form 
$\partial_t u = \Delta u + a(x)u^\mathsf{p} + b(x) u^\mathsf{q} + c(x) u$ with 
$\mathsf{p}=(n+2)/(n-2)$ and $\mathsf{q}=(3n-2)/(n-2)$ while when $n=2$ the same evolutionary equation takes 
the form $\partial_t u = \Delta u + a(x) e^{2u} + b(x) e^{-2u} + c(x)$. In the setting of smooth metric measure spaces 
a generalisation of the Einstein-scalar field Lichnerowicz equation with space-time dependent coefficients can be 
described as:  
\begin{equation}
\square_q u = [\partial_t -q(x,t) - \Delta_f] u = a(x,t) u^\mathsf{p} + b(x,t) u^\mathsf{q} + c(x,t) u \log u, 
\end{equation}
and 
\begin{equation}
\square u = (\partial_t - \Delta_ f) u = a(x,t) e^{2u} + b(x,t) e^{-2u} + c(x,t). 
\end{equation}
For gradient estimates, Harnack inequalities, Liouville type theorems and other related results in this direction 
see \cite{Dung, Song, Taheri-GE-1, Taheri-GE-2, [Wu18]} and the references therein.

Moving on to the evolving case the time dependence of the metric-potential pair adds further 
complications and technical details as far as gradient estimates are concerned. Here the case of the weighted 
heat equation under the Perelman-Ricci flow, generalising in turn, the heat equation under the Ricci 
flow to the setting to smooth metric measure spaces, given by the system  
\begin{align} \label{eqRicciFlow} 
\begin{cases}
\square u(x,t) = \left[ \dfrac{\partial}{\partial t} - \Delta_f \right] u(x,t) =0, \\
\dfrac{\partial g}{\partial t}(x,t) = - 2 {\mathscr Ric}_f (g)(x,t), \\ 
\end{cases}
\end{align}
has been considered by many authors (see in particular 
\cite{[Ba1], CaoChow, Chow, [Ha93], Ham, MReto, [Pe02], SunJ, [QZh]} 
and the references therein).

The system \eqref{eq11} can be seen as a generalisation of \eqref{eqRicciFlow} in two substantial ways. 
Firstly, the weighted linear heat equation is replaced by its nonlinear counterpart where the nonlinearity 
takes a considerably general formulation. Secondly, the Perelman-Ricci flow gives way to the 
$(\mathsf{k},m)$-{\it super} Perelman-Ricci flow which is again a substantial and far-reaching 
generalisation. For example, in the static case (with $\partial_t g \equiv 0$, $\partial_t f \equiv 0$), 
the latter includes spaces with curvature lower bound ${\mathscr Ric}_f(g) \ge -\mathsf{k} g$ 
(see the discussion in the next section) whereas the counterpart of the former [{\it cf}. \eqref{eqRicciFlow}] 
in the static case includes only ${\mathscr Ric}_f$-flat spaces which is a much narrower subset. 
In passing we point out that ${\mathscr Ric}$-flat spaces are {\it particular} classes of Einstein 
manifolds, verifying ${\mathscr Ric}(g) = \lambda g$ with $\lambda=0$, 
and ${\mathscr Ric}_f^m$-flat spaces, naturally extending the notion to 
the smooth metric measure space context, are particular classes of 
$m$-quasi Einstein manifolds, verifying ${\mathscr Ric}_f^m(g)=\lambda g$, 
or gradient Ricci solitons, verifying ${\mathscr Ric}_f(g) = \lambda g$, 
as indicated earlier in \eqref{eq1.5}, each in the special case with 
$\lambda=0$ ({\it see} \cite{[LLi2], [LLi4], [QZh]}).

For the sake of future reference a $(\mathsf{k},m)$-{\it super} Perelman-Ricci flow specifically refers to a 
complete smooth solution $(g,f)$ to the flow inequality (with $n \le m < \infty$):
\begin{align} \label{eq11b} 
\begin{cases}
\dfrac{\partial g}{\partial t}(x,t) + 2 {\mathscr Ric}_f^m(g)(x,t) \ge - 2 \mathsf{k}g(x,t), \\ 
{\mathscr Ric}_f^m(g)(x,t) = {\mathscr Ric} (g)(x,t) + {\rm Hess} (f) (x,t) 
- \dfrac{\nabla f \otimes \nabla f}{m-n} (x,t).   
\end{cases}
\end{align}

In the event $m=\infty$ and with ${\mathscr Ric}_f(g)$ as described by \eqref{Ricf def eq} 
the above system should be interpreted explicitly as 
\begin{align} \label{eq11c} 
\begin{cases}
\dfrac{\partial g}{\partial t}(x,t) + 2 {\mathscr Ric}_f(g)(x,t) \ge - 2 \mathsf{k}g(x,t), \\ 
{\mathscr Ric}_f(g)(x,t) = {\mathscr Ric} (g)(x,t) + {\rm Hess} (f) (x,t).    
\end{cases}
\end{align}
Thus hereafter when referring to the flow inequality in \eqref{eq11} (the symmetric tensor 
inequality on the second line) one of the above is intended depending on whether $m$ is finite or not.

\qquad \\
{\bf Plan of the paper.} Let us finish-off this long introduction by briefly describing the layout and plan of the paper. 
In Section \ref{sec2} we present the statements 
of the main results. Section \ref{sec3} is devoted to the proof of the local Soulpet-Zhang 
estimate in Theorem \ref{thm1} followed by the proof of the elliptic Harnack inequality in Theorem \ref{cor Harnack}. 
In Section \ref{sec5} we present the proof of the local Hamilton-type 
estimates in Theorem \ref{thm18} and in 
Section \ref{sec6}, which is the core and most involved part of the paper, we establish 
the local parabolic Li-Yau type estimate in Theorem \ref{thm28}. 
In Section \ref{sec7} we present the proof of the parabolic Harnack inequality in 
Theorem \ref{thm38} and in Section \ref{sec8} we give the proof of the two Liouville 
results in Theorem \ref{thm46} and Theorem \ref{thm48}. 
Finally in Section \ref{sec9} we present the proof of the 
global bound in Theorem \ref{thm25} and discuss its consequences.

\qquad \\
{\bf Notation.} For $X \in \mathbb{R}$ we write $X_+=\max(X, 0)$ and $X_-=\min(X, 0)$. 
Therefore $X=X_+ + X_-$ with $X_+ \ge 0$ and $X_- \le 0$. 
Fixing a reference point $x_0 \in M$ we denote by $d=d(x,x_0, t)$ the Riemannian distance between $x$ 
and $x_0$ on $M$ with respect to the metric $g=g(t)$. We write $\varrho=\varrho(x,x_0,t)$ for the geodesic 
radial variable measuring the distance between $x$ and $x_0$ at time $t>0$. For $R>0$, $T>0$ we define 
the space-time cylinder 
\begin{equation} \label{space-time-cylinder}
Q_{R,T}(x_0) \equiv \{ (x, t) : d(x, x_0, t) \le R \mbox{ and } 0 \le t \le T \} \subset M \times [0, T], 
\end{equation}
and for fixed $0<t \le T$ we denote by $B_r(x_0) \subset M$ the geodesic ball of radius $r>0$ centred at 
$x_0$ with respect to $g=g(t)$. When the choice of the reference point $x_0$ is clear from the context we 
often abbreviate and simply write $d(x, t)$, $\varrho(x,t)$ or $B_r$, $Q_{R,T}$ respectively.

We typically denote partial derivatives by subscripts unless otherwise specified. In particular, for the nonlinear 
function $\Sigma=\Sigma(t,x,u)$ in \eqref{eq11} we frequently use $\Sigma_x$, $\Sigma_u$, $\Sigma_{xx}$, 
$\Sigma_{xu}$ and $\Sigma_{uu}$ for the respective partial derivatives of first and second orders in the spatial 
variables $x$ and $u$. [Note that in local coordinates we have $\Sigma_x=(\Sigma_{x_1}, \dots, \Sigma_{x_n})$.]
We also write $\Sigma^x: x \mapsto \Sigma(t,x, u)$ for the function resulting from freezing the variables $(t,u)$ 
and viewing $\Sigma$ as a function of $x$ only. Thus we speak of $\nabla \Sigma^x$, $\Delta \Sigma^x$, 
$\Delta_f \Sigma^x$ and so on.

For a bounded function $u=u(x,t)$ on $Q_{R,T}$ we write $\overline u=\overline u (R,T)$ and 
$\underline u=\underline u (R, T)$ for the supremum and infimum of $u$ on $Q_{R, T}$ 
respectively.  We also introduce the set
\begin{equation}
\Theta_{R, T} \equiv \{(t,x,u) : (x,t) \in Q_{R, T} \mbox{ and } \underline u \le u \le \overline u \} 
\subset [0, T] \times M \times \mathbb{R}.
\end{equation} 
It is evident that for any $F=F(t,x,u)$ we have the inequalities 
\begin{equation}
\sup_{Q_{R, T}} F(t,x,u(x,t)) \le \sup_{\Theta_{R,T}} F(t,x,u), 
\end{equation}
\begin{equation}
\inf_{Q_{R, T}} F(t,x,u(x,t)) \ge \inf_{\Theta_{R,T}} F(t,x,u).  
\end{equation}

Note that whilst the quantities on the left depend explicitly on the function $u$, the quantities on the right 
depend only on the upper and lower bounds $\overline u$, $\underline u$ of $u$. Depending on context 
and need both these bounds will be utilised as appropriate in future.

\section{Statement of the main results} \label{sec2}

In this section we present the main results of the paper along with some accompanying discussion. The complete 
proofs and further details are delegated to the subsequent sections. For the convenience of the reader and 
ease of reference we have grouped this section into four subsections based on the nature of the estimates 
and results involved.

\subsection{A local and a global Souplet-Zhang type gradient estimate for \eqref{eq11}}

The first result is a local elliptic gradient estimate of Souplet-Zhang type for positive smooth bounded solutions 
to \eqref{eq11}. We emphasise that here the metric-potential pair $(g,f)$ is assumed to be time dependent and 
a complete smooth solution to the flow inequality \eqref{eq11c} for a suitable $\mathsf{k} \ge 0$. Note that as a 
result of this evolution the usual differential operators $\nabla$, ${\rm div}$, $\Delta$ and 
$\Delta_f = \Delta - \langle \nabla f, \nabla \rangle$ are all time dependent.

For the sake of this local estimate we pick a reference point $x_0 \in M$ and restrict to the compact set 
$Q_{R,T}=Q_{R,T}(x_0)$ where $R \ge 2$ and $T>0$ are fixed. 
The estimate makes use of the upper bound on the solution $u$, the lower bound on 
the generalised Bakry-\'Emery Ricci curvature tensor ${\mathscr Ric}_f(g) \ge -(n-1) k_1 g$ 
and the lower bound $\partial_t g \ge -2k_2 g$, with $k_1, k_2 \ge 0$, 
all within the set $Q_{R,T}$. Two important quantities appearing in the formulation of the estimate that directly 
link to the nonlinearity $\Sigma$ and the solution $u$ are respectively $\mathsf{R}_\Sigma = \mathsf{R}_\Sigma(u)$, 
$\mathsf{P}_\Sigma = \mathsf{P}_\Sigma(u)$ [see \eqref{eq13}-\eqref{eq13-R}]. 
These play a key role in the subsequent bounds, Harnack inequalities and Liouville-type results on solutions. 
Note also that as $u>0$ and $Q_{R,T}$ is compact, $u$ is bounded away from zero 
and from above hence these quantities are finite.

\begin{theorem} \label{thm1}
Let $(M, g, d\mu)$ be a complete smooth metric measure space with $d\mu=e^{-f} dv_g$ and suppose that the 
metric-potential pair $(g, f)$ is time dependent and of class $\mathscr{C}^2$. Assume the bounds
\begin{equation}
{\mathscr Ric}_f (g) \ge -(n-1) k_1 g, \qquad \partial_t g \ge -2 k_2 g, 
\end{equation} 
in the compact set $Q_{R,T}$ for some $k_1, k_2 \ge 0$. Let $u$ be a positive solution to \eqref{eq11} 
with $0<u \le D$ in $Q_{R,T}$. Then for all $(x,t)$ in $Q_{R/2,T}$ with $t>0$ we have:
\begin{align}\label{eq13}
\frac{|\nabla u|}{u} 
\le C \left\{ \frac{1}{R} + \sqrt{\frac{[\gamma_{\Delta_f}]_+}{R}} + \frac{1}{\sqrt{t}} + \sqrt{k} 
+ \sup_{Q_{R,T}} \left( \mathsf{N}_q + \mathsf{R}^{1/2}_\Sigma(u) 
+ \mathsf{P}^{1/3}_\Sigma(u) \right) \right\} \left(1 - \log \frac{u}{D} \right), 
\end{align}
where $C>0$ depends only on $n$, $\mathsf{N}_q = q_+^{1/2} + |\nabla q|^{1/3}$ with  
$q_+ = q_+(x,t) = [q(x,t)]_+$ and 
\begin{align}\label{eq13-R}
\mathsf{R}_\Sigma(u) = \left[ \frac{\Sigma_u(t,x,u)}{1-\log(u/D)}   + \frac{\log (u/D) \Sigma(t,x,u)}{u[1-\log(u/D)]^2} \right]_+, \quad 
\mathsf{P}_\Sigma(u) = \frac{|\Sigma_x(t,x,u)|}{u[1-\log (u/D)]^2}.
\end{align}
Additionally, $k= \sqrt{k_1^2 + k_2^2}$ and $[\gamma_{\Delta_f}]_+=\max(\gamma_{\Delta_f}, 0)$ where  
\begin{equation} \label{sigma def}
\gamma_{\Delta_f} = \max_{(x, t)} \{ \Delta_f \varrho (x,t) : d (x,x_0, t) =1 \mbox{ and } 0 \le t \le T \}.
\end{equation}
\end{theorem}

Let us make a few useful comments on the theorem, its assumptions and proof. 
First, ${\mathscr Ric}_f(g) \ge - (n-1)k_1 g$ and $\partial_t g \ge - 2k_2 g$ in $Q_{R,T}$ give 
\eqref{eq11c} with $\mathsf{k} = (n-1)k_1 + k_2 \ge 0$. In the proof, the lower 
bound on ${\mathscr Ric}_f(g)$ is utilised in the application of the generalised Laplacian comparison theorem and 
the bound $\partial_t g \ge - 2k_2 g$ is used for controlling the time derivative of the geodesic distance function: 
$\partial_t d= \partial [d(x,x_0,t)]/\partial t$. Both these estimates arise in the localisation stage in the later part 
of the proof.

Second, in the static case we have $\partial_t g \equiv 0$, $\partial_t f \equiv 0$ and ${\mathscr Ric}_f(g) \ge -(n-1) k_1 g$ 
and so we can set $k_2=0$ and $k=k_1$. [This case is certainly a solution to the flow inequality in \eqref{eq11} with 
$m=\infty$ and $\mathsf{k}=(n-1)k_1$.] Thus here $\gamma_{\Delta_f} = \max \{ \Delta_f \varrho (x) : d(x,x_0)=1\}$ and 
evidently $\partial_t d= \partial [d(x,x_0,t)]/\partial t \equiv 0$. This way Theorem \ref{thm1} can also be seen 
as giving local gradient estimates for positive bounded solutions to \eqref{eq11} on non-evolving (static) smooth metric measure 
spaces which is of course of independent interest. Finally, note that by virtue of the bound $0<u \le D$ we have the inequality 
$0<1/[1-\log(u/D)]\le1$ and so in particular 
\begin{equation} \label{compare}
0 \le \mathsf{R}_\Sigma(u) \le \left[ \frac{u\Sigma_u(t,x,u) + \log (u/D) [\Sigma(t,x,u) - u \Sigma_u(t,x,u)]}{u [1- \log (u/D)]} \right]_+, 
\end{equation} 
and $0 \le \mathsf{P}_\Sigma (u) \le  |\Sigma_x(t,x,u)|/(u[1-\log(u/D)]) \le |\Sigma_x(t,x,u)|/u$. (It is instructive to compare the 
bound \eqref{compare} with those in \cite{Taheri-GE-1, Taheri-GE-2}.)

The local estimate above has a global in space counterpart subject to the prescribed bounds in the theorem being global in space. 
The proof follows by passing to the limit $R \to \infty$ in \eqref{eq13} and taking into account the vanishing of the terms involving 
$R$ by virtue of the bounds being global and the relevant constants being independent of $R$. The precise formulation of this is 
given in the following theorem.

\begin{theorem} \label{thm1-global}
Let $(M, g, d\mu)$ be a complete smooth metric measure space with $d\mu=e^{-f} dv_g$ and suppose that the 
metric-potential pair $(g, f)$ is time dependent and of class $\mathscr{C}^2$. Assume ${\mathscr Ric}_f (g) \ge -(n-1) k_1 g$ 
and $\partial_t g \ge -2 k_2 g$ in $M \times [0, T]$ for $k_1, k_2 \ge 0$. If $u$ is a positive 
solution to \eqref{eq11} with $0<u \le D$, then for all $x \in M$ and $0<t \le T$ we have:
\begin{align}\label{eq13-global}
\frac{|\nabla u|}{u} 
\le C \left\{ \frac{1}{\sqrt{t}} + \sqrt{k} 
+ \sup_{M \times [0, T]} \left( \mathsf{N}_q + \mathsf{R}^{1/2}_\Sigma(u) 
+ \mathsf{P}^{1/3}_\Sigma(u) \right) \right\} \left(1 - \log \frac{u}{D} \right).  
\end{align}
The quantities $\mathsf{N}_q$, $\mathsf{R}_\Sigma(u)$ and $\mathsf{P}_\Sigma(u)$ in \eqref{eq13-global} 
are the same as those in Theorem $\ref{thm1}$.
\end{theorem}

One of the useful consequences of the estimates above is the following elliptic Harnack inequality for bounded positive 
solutions to equation \eqref{eq11}. Later on we will also prove a parabolic counterpart for this inequality from another 
type of estimate on solutions ({\it see} Theorem \ref{thm38} below and Section \ref{sec7} for the proof). Note that here 
the solution $u$ is compared at two different spatial points $x_1$, $x_2$ but the same time $t>0$.

\begin{theorem} \label{cor Harnack}
Under the assumptions of Theorem $\ref{thm1}$ for all $(x_1,t)$, $(x_2, t)$ in $Q_{R/2,T}$ with $t>0$ we have:
\begin{equation} \label{eq2.5}
 u(x_1, t) \le (eD)^{1-\gamma} [u(x_2, t)]^\gamma, 
\end{equation}
where the exponent $\gamma$ in \eqref{eq2.5} is explicitly given $($with $d=d(x_1, x_2, t)$ below$)$ by
\begin{equation} \label{Exp-gamma-2.7}
\gamma = {\rm exp} \left[ - C d 
\left( \frac{1}{R}  
+ \sqrt k + \frac{1}{\sqrt{t}}  
+ \sup_{Q_{R,T}} \left[ \mathsf{N}_q + \mathsf{R}^{1/2}_\Sigma(u) + \mathsf{P}^{1/3}_\Sigma(u) \right] 
+ \sqrt{ \frac{[\gamma_{\Delta_f}]_+}{R}}
\right) \right].
\end{equation}
Moreover, subject to the global bounds in Theorem $\ref{thm1-global}$, for all $x_1$, $x_2$ in $M$ and $0<t \le T$ 
we have the same \eqref{eq2.5}-\eqref{Exp-gamma-2.7} with $M \times [0, T]$ replacing $Q_{R,T}$ in \eqref{Exp-gamma-2.7}. 
\end{theorem}

\subsection{A local and a global Hamilton-type gradient estimate for \eqref{eq11}}

We now present the second gradient estimate of elliptic type for the positive smooth solutions 
to \eqref{eq11}. Here again the metric-potential pair $(g,f)$ is time dependent and a complete smooth solution 
to the flow inequality \eqref{eq11c} and the estimate makes use of the bounds on the solution $u$, the 
lower bound on the generalised Bakry-\'Emery Ricci curvature tensor ${\mathscr Ric}_f(g) \ge -(n-1) k_1 g$ 
and the lower bound $\partial_t g \ge -2k_2 g$, with $k_1, k_2 \ge 0$, all within the compact set $Q_{R,T}$. 
The setting and notation here is mostly similar to those in the previous 
theorem however the proof is based on the use of different objects and tools.

\begin{theorem} \label{thm18}
Let $(M, g, d\mu)$ be a complete smooth metric measure space with $d\mu=e^{-f} dv_g$ and suppose that the 
metric-potential pair $(g, f)$ is time dependent and of class $\mathscr{C}^2$. Assume the bounds 
\begin{equation}
{\mathscr Ric}_f (g) \ge -(n-1) k_1 g, \qquad \partial_t g \ge -2 k_2 g,
\end{equation} 
in the compact set $Q_{R,T}$ for some $k_1, k_2 \ge 0$. Let $u$ be a positive 
solution to \eqref{eq11} in $Q_{R,T}$. Then for all $(x,t)$ in $Q_{R/2,T}$ with $t>0$ we have:
\begin{align} \label{eq1.12}
\frac{|\nabla u|}{\sqrt u} \le C \left\{ \frac{1}{R} + \sqrt{\frac{[\gamma_{\Delta_f}]_+}{R}} + \frac{1}{\sqrt t} + \sqrt k 
+ \sup_{Q_{R, T}} \left( \mathsf{N}_q + \mathsf{T}_\Sigma^{1/2} (u)+ \mathsf{S}_\Sigma^{1/3} (u) \right) \right\} 
\Big( \sup_{Q_{R, T}} \sqrt u \Big),  
\end{align}
where $C>0$ depends only on $n$, $\mathsf{N}_q = q_+^{1/2} + |\nabla q|^{1/3}$, 
$\gamma_{\Delta_f}$ is as in \eqref{sigma def}, $k = \sqrt{k_1^2 + k_2^2}$ and  
\begin{equation} \label{eq2.8}
\mathsf{T}_\Sigma (u) = \left[ \frac{2u \Sigma_u(t,x,u)-\Sigma(t,x,u)}{u} \right]_+, \qquad \mathsf{S}_\Sigma(u) = \frac{|\Sigma_x (t,x,u)|}{u}.  
\end{equation}
\end{theorem}

It is instructive to compare the nonlinear quantities $\mathsf{R}_\Sigma(u)$, $\mathsf{P}_\Sigma(u)$ in \eqref{eq13-R} 
of Theorem \ref{thm1} with the corresponding ones $\mathsf{T}_\Sigma(u)$, $\mathsf{S}_\Sigma(u)$ in 
\eqref{eq2.8} of Theorem \ref{thm18} as appearing on the right-hand sides of the formulations of the estimates for the gradient of the solution $u$. 
More on this and its implications will be said later when discussing Liouville-type results, Harnack inequalities and other applications of the estimates.

Again, the local estimate above has a global counterpart, when the asserted bounds in the theorem are global. 
The proof follows by passing to the limit $R \to \infty$ in \eqref{eq1.12} taking into account the vanishing of the 
terms involving $R$ resulting from the bounds being global and the relevant constants being independent of $R$. 
This is the content of the following theorem.

\begin{theorem} \label{thm18-global}
Let $(M, g, d\mu)$ be a complete smooth metric measure space with $d\mu=e^{-f} dv_g$ and suppose that the 
metric-potential pair $(g, f)$ is time dependent and of class $\mathscr{C}^2$. Assume ${\mathscr Ric}_f (g) \ge -(n-1) k_1 g$ 
and $\partial_t g \ge -2 k_2 g$ in $M \times [0, T]$ for $k_1, k_2 \ge 0$. If $u$ is a positive solution to \eqref{eq11}, 
then for all $x \in M$ and $0<t \le T$ we have:
\begin{align} \label{eq1.12-global}
\frac{|\nabla u|}{\sqrt u} \le C \left\{\frac{1}{\sqrt t} + \sqrt k 
+ \sup_{M \times [0, T]} \left( \mathsf{N}_q + \mathsf{T}_\Sigma^{1/2} (u)+ \mathsf{S}_\Sigma^{1/3} (u) \right) \right\} 
\Big( \sup_{M \times [0, T]} \sqrt u \Big). 
\end{align}
The quantities $\mathsf{N}_q$, $\mathsf{T}_\Sigma (u)$ and $\mathsf{S}_\Sigma (u)$ in \eqref{eq1.12-global} 
are the same as those in Theorem $\ref{thm18}$. 
\end{theorem}

\subsection{A local and a global parabolic differential Harnack inequality for \eqref{eq11}}

We now move on to a different type of estimate to the elliptic ones established above, namely, 
a Li-Yau type estimate (also known as a differential Harnack inequality). To this end we pick 
$x_0 \in M$ and set $Q_{2R,T}=Q_{2R,T}(x_0)$ where $R>0$, $T>0$. 
The estimate makes use of the lower bound on ${\mathscr Ric}^m_f(g)$ (with $n \le m < \infty$), along with the bounds 
on $u$, $\partial_t g$, $\nabla f$, $\nabla \partial_t g$ and $\nabla \partial_t f$ all 
in $Q_{2R,T}$. 
Applications to Harnack inequalities and Liouville-type results for solutions of \eqref{eq11} are given later on.

\begin{theorem} \label{thm28}
Let $(M, g, d\mu)$ be a complete smooth metric measure space with $d\mu=e^{-f} dv_g$ and suppose that the 
metric-potential pair $(g, f)$ is time dependent and of class $\mathscr{C}^2$. 
Assume that ${\mathscr Ric}_f^m(g) \ge -(m-1) k_1 g$ and that the following bounds 
\begin{align} \label{k-bounds-one}
-2\underline k_2 g \le \partial_t g  \leq 2 \overline k_2 g, \qquad |\nabla \partial_t g| \le 2k_3, 
\end{align}
\begin{align} \label{k-bounds-two}
|\nabla f| \le \ell_1, \qquad |\nabla \partial_t f| \le \ell_2, 
\end{align} 
hold in $Q_{2R,T}$ for suitable constants $k_1, \underline k_2, \overline k_2, k_3 \ge 0$ and $\ell_1, \ell_2 \ge 0$. 
Let $u=u(x,t)$ be a positive solution to \eqref{eq11} in $Q_{2R,T}$. Then for every $\lambda>1$, 
$\varepsilon \in (0, 1)$ and all $(x,t)$ in $Q_{R,T}$ with $t>0$ we have 
\begin{align} \label{1.26}
\frac{|\nabla u|^2}{\lambda u^2} - \frac{\partial_t u}{u} + q + \frac{\Sigma(t,x,u)}{u}  
\le&~m \lambda \left[ \frac{1}{t} + c_1 \underline k_2 + \gamma^\Sigma_1 \right] \nonumber \\
&+ \frac{m \lambda}{R^2} \left[ \frac{m c_1 ^2 \lambda ^2}{2(\lambda -1)}+c_2 +(m-1) c_1(1+R \sqrt{k_1})+2c_1^2 \right]  \nonumber \\
&+ \sqrt m \bigg\{ \frac{m \lambda^2 \mathsf{A}^2}{4(1-\varepsilon)(\lambda-1)^2} 
+ \frac{3}{4} \left[ \frac{m \lambda^2 \mathsf{B}^4}{4 \varepsilon(\lambda-1)^2} \right]^{1/3} 
\nonumber \\
&+ \lambda ^2 n (\underline k_2 + \overline k_2)^2 + 2 \lambda^2 n k_3 
+ \lambda (\gamma^\Sigma_2 + \gamma^{qu}_2) \bigg\}^{1/2}.  
\end{align}
The quantities and constants appearing on the right-hand side of the inequality \eqref{1.26} are given respectively by 
\begin{align} \label{eq2.13}
\mathsf{A} =&~2[(m-1)k_1 +(\lambda -1) \overline k_2 + k_3] \nonumber \\
&- \inf_{\Theta_{2R, T}} \left\{ \frac{1}{u} [ \Sigma(t,x,u) - u\Sigma_u(t,x,u) + \lambda u^2 \Sigma_{uu}(t,x,u) ]_- \right\},  
\end{align}
and
\begin{equation} \label{eq2.15}
\mathsf{B} = \lambda \ell_2 + 2 \lambda \underline k_2 \ell_1 
+ \sup_{\Theta_{2R, T}}  \left\{ \frac{2}{u} |\Sigma_x(t,x,u) -\lambda u \Sigma_{xu}(t,x,u)| + 2( \lambda -1) |\nabla q| \right\}.  
\end{equation}
Moreover we have 
\begin{equation} \label{eq2.16}
\gamma^\Sigma_1 = \sup_{\Theta_{2R, T}} \left\{ \frac{1}{u} [ u \Sigma_u (t,x,u) - \Sigma(t,x,u) ]_{+} \right\}, 
\end{equation}
and similarly 
\begin{equation} \label{eq2.17}
\gamma^\Sigma_2 = - \inf_{\Theta_{2R, T}} \left[ \frac{1}{u} \Delta_f \Sigma^x (t,x,u) \right]_-, \qquad 
\gamma^{qu}_2 = - \inf_{\Theta_{2R, T}} \left[ \Delta_f q \right]_-.
\end{equation} 
Finally $\Theta_{2R, T} = \{(t,x,u) : (x,t) \in Q_{2R, T}, \, \underline u \le u \le \overline u \} \subset [0, T] \times M \times (0, \infty)$, 
where $\overline u$, $\underline u$ denote the maximum and minimum of $u$ on the compact space-time cylinder $Q_{2R, T}$. 
\end{theorem}

Note that unlike in Theorem \ref{thm1} and Theorem \ref{thm18}, here, 
the curvature lower bound is imposed on the generalised Bakry-\'Emery Ricci tensor ${\mathscr Ric}_f^m(g)$ 
with $n \le m <\infty$ and not on ${\mathscr Ric}_f(g)$. Evidently, a lower bound on ${\mathscr Ric}^m_f(g)$ 
is a stronger condition than one on ${\mathscr Ric}_f(g)$ as the former implies the latter but not {\it vice versa} 
[see \eqref{Ricci-m-f-intro}-\eqref{Ricf def eq}]. Also as is seen from \eqref{Bochner}, \eqref{Bochner-m-inequality} 
the bound ${\mathscr Ric}^m_f(g) \ge -\mathsf{k} g$ leads to the curvature-dimension condition 
${\rm CD}(-\mathsf{k},m)$ [whilst ${\mathscr Ric}_f(g) \ge -\mathsf{k} g$ leads to ${\rm CD}(-\mathsf{k}, \infty)$] 
({\it cf.} \cite{BE, BD, Bak}, \cite{[Wu18], [WuW14]}). From a purely technical view point these two conditions 
are of fundamentally different strengths and nature ({\it see} \cite{Bak} for more on their implications, in particular, on functional 
and other types of geometric inequalities associated with their respective diffusion operators and semigroups).

The global counterpart of the above local estimate can be obtained by imposing appropriate global bounds in the assumptions 
as is formulated in the following theorem.

\begin{theorem} \label{thm28-global}
Let $(M, g, d\mu)$ be a complete smooth metric measure space with $d\mu=e^{-f} dv_g$ and assume that the 
metric-potential pair $(g, f)$ is time dependent and of class $\mathscr{C}^2$. Assume 
${\mathscr Ric}^m_f (g) \ge -(m-1) k_1 g$, \eqref{k-bounds-one} and \eqref{k-bounds-two} hold globally 
in $M \times [0, T]$. 
Let $u=u(x,t)$ be a positive solution to \eqref{eq11}. Then for every 
$\lambda>1$, $\varepsilon \in (0, 1)$ and all $x \in M$, $0<t \le T$ we have 
\begin{align} \label{1.26-global}
\frac{|\nabla u|^2}{\lambda u^2} - \frac{\partial_t u}{u} + q + \frac{\Sigma(t,x,u)}{u} 
\le&~m \lambda \left[ \frac{1}{t} + c_1 \underline k_2 + \gamma^\Sigma_1 \right] \nonumber \\
&+ \sqrt m \bigg\{ \frac{m \lambda^2 \mathsf{A}^2}{4(1-\varepsilon)(\lambda-1)^2} 
+ \frac{3}{4} \left[ \frac{m \lambda^2 \mathsf{B}^4}{4 \varepsilon(\lambda-1)^2} \right]^{1/3} \nonumber \\
&+ \lambda ^2 n (\underline k_2 + \overline k_2)^2 + 2 \lambda^2 n k_3 
+ \lambda (\gamma^\Sigma_2 + \gamma^{qu}_2) \bigg\}^{1/2}. 
\end{align}
Here $\mathsf{A}$ is as in \eqref{eq2.13}, $\mathsf{B}$ is as in \eqref{eq2.15} and the quantities $\gamma^\Sigma_1$, 
$\gamma^\Sigma_2$ and $\gamma_2^{qu}$ are as in \eqref{eq2.16} and \eqref{eq2.17} with the supremum and infimums taken over 
$M \times [0, T]$ respectively. 
\end{theorem}

\begin{theorem} \label{thm38}
Under the assumptions of Theorem $\ref{thm28}$, if $u$ is a positive solution to \eqref{eq11}, 
then for every $(x_1, t_1)$, $(x_2, t_2)$ in $Q_{R, T}$ with $t_2>t_1$ and $\lambda >1$ we have 
\begin{align}
u(x_2,t_2)\geq u(x_1,t_1)\left(\frac{t_2}{t_1}\right)^{-m \lambda}e^{- \lambda L(x_1,x_2, t_2-t_1)} e ^ {(t_2-t_1) S}.  
\end{align}
Here $S$ is a constant depending only on the bounds given in Theorem $\ref{thm28}$ and is explicitly given 
by \eqref{S-eq7.1}. Furthermore $L$ is given by 
\begin{align}
L(x_1,x_2, t_2-t_1) = \inf_{\gamma \in \Gamma}  \left[ \frac{1}{4(t_2-t_1)} \int_{0}^{1} |\dot \gamma(t)|_{g(t)}^2\,dt \right], 
\end{align} 
where $\Gamma$ is the set of all curves $\gamma \in \mathscr{C}^1( [t_1,t_2]; M)$ lying entirely in $Q_{R, T}$ 
with $\gamma(t_1) = x_1$ and $\gamma(t_2) = x_2$. 
If the bounds are global as in Theorem $\ref{thm28-global}$ then the above estimate is global. 
\end{theorem}

\subsection{Liouville-type results, global bounds and elliptic Harnack inequalities}

We now move on to presenting some more applications of the results formulated above. Towards this end 
we begin with two independent global constancy and Liouville-type results for the elliptic counterpart 
of equation \eqref{eq11}. The first result follows from the local elliptic gradient estimate in Theorem \ref{thm18}
and the second result puts to use the local parabolic gradient estimate in Theorem \ref{thm28}. 
Some closely related applications, in particular, to parabolic Liouville-type results and ancient 
solutions of \eqref{eq11} that are equally of great interest are discussed in our forthcoming paper.

\begin{theorem} \label{thm46}
Let $(M, g, d\mu)$ be a complete smooth metric measure space with $d\mu=e^{-f} dv_g$ and time independent 
metric and potential of class $\mathscr{C}^2$. Assume ${\mathscr Ric}_f(g) \ge 0$ in $M$. Let $u=u(x)$ be 
a positive bounded solution to the equation 
\begin{equation}
\Delta_ f u + \Sigma(u) = 0, 
\end{equation}
such that $\Sigma(u) - 2u \Sigma_u(u) \ge 0$ everywhere in $M$. Then $u$ must be a constant. 
\end{theorem}

\begin{theorem} \label{thm48}
Let $(M, g, d\mu)$ be a complete smooth metric measure space with $d\mu=e^{-f} dv_g$ and time independent 
metric and potential of class $\mathscr{C}^2$. Assume ${\mathscr Ric}_f^m(g) \ge 0$ in $M$. Let $u=u(x)$ be 
a positive solution to 
\begin{equation}
\Delta_ f u + \Sigma(u) = 0.
\end{equation}
Then for every $\lambda>1$, $\varepsilon \in (0, 1)$ and all $x \in M$ we have 
\begin{align} \label{eqL2.26}
\frac{|\nabla u|^2}{\lambda u^2} + \frac{\Sigma(u)}{u} 
\le&~ m \lambda 
\sup_{\Theta} \left\{ \left[\frac{-[\Sigma(u) - u \Sigma_u (u)]}{u}\right]_{+} \right\} \nonumber \\
&+ \frac{m \lambda/\sqrt{1-\varepsilon}}{2(\lambda-1)} 
\sup_{\Theta} \left\{ \left[ \frac{- [\Sigma(u) - u\Sigma_u(u) + \lambda u^2 \Sigma_{uu}(u)]}{u} \right]_+ \right\}.  
\end{align}
In particular if along the solution $u$ we have $\Sigma(u) \ge 0$, $\Sigma(u) - u \Sigma_u(u) \ge 0$ and for some 
$\lambda>1$ 
\begin{equation}
\Sigma(u)-u\Sigma_u(u)+\lambda u^2 \Sigma_{uu}(u) \ge 0, 
\end{equation}
everywhere on $M$ then $u$ must be a constant. 
\end{theorem}

The final result we present here is a Hamilton-type global bound and a corresponding global Harnack 
interpolation inequality for positive solutions to \eqref{eq11} under the flow inequality \eqref{eq11c}. 
Here we assume $M$ is closed.

\begin{theorem} \label{thm25}
Let $u$ be a positive solution to \eqref{eq11} with $q=0$ satisfying $0<u \le D$ in 
$M \times [0,T]$ and assume that the metric-potential pair $(g,f)$ evolves under the $\mathsf{k}$-super Perelman-Ricci 
flow \eqref{eq11c} with $\mathsf{k} \ge 0$. Assume that $\Sigma(u) \ge 0$ and $\Sigma(u)-2u\Sigma_u(u)\ge 0$ along the 
solution $u$. Then 
\begin{align} \label{eq15}
t |\nabla \log u|^2 \le (1+2\mathsf{k}t) [1+ \log (D/u)], \qquad 0<t<T. 
\end{align}
As a result there holds the following Harnack-interpolation inequality: For any $s>0$, $x_1, x_2 \in M$ and $0<t \le T$, 
we have 
\begin{align}\label{eq16}
u(x_1,t) \le (eD)^{s/(1+s)} \exp \left( d^2(x_1,x_2,t) \frac{1+2kt}{4st} \right) [u(x_2, t)]^{1/(1+s)}.
\end{align}
\end{theorem}

\section{Proof of the Souplet-Zhang estimate in Theorem \ref{thm1}}\label{sec3}

\subsection{Some intermediate parabolic lemmas (${\bf I}$)}

Before proceeding onto presenting the proof of Theorem \ref{thm1} we need some intermediate results 
and lemmas. In this subsection we prove a parabolic differential inequality utilised in the proof of the 
theorem. This is inequality \eqref{eq2112} in Lemma \ref{lem2112} where $u$ is a positive solution 
to \eqref{eq11} with the metric-potential pair $(g,f)$ being a complete solution to the $\mathsf{k}$-super 
Perelman-Ricci flow \eqref{eq11c}.  We begin with the following parabolic differential {\it identity} first.

\begin{lemma} \label{lem20}
Let $u$ be a positive bounded solution to the equation $(\partial_t - \Delta_f) u = \Sigma (t,x,u)$ with 
$0<u \le D$. Then the function $h = \log (u/D)$ satisfies the equation 
\begin{equation} \label{h evolution eq}
\square h = (\partial_t - \Delta_f) h = |\nabla h|^2 + D^{-1} e^{-h} \Sigma(t,x,D e^h).
\end{equation} 
\end{lemma}

\begin{proof} 
This is an easy calculation and the proof is left to the reader. 
\end{proof}

\begin{lemma} \label{lem21}
Suppose $g=g(t)$, $f=f(t)$ are of class $\mathscr{C}^2$ and let $u$ be a positive bounded solution to the equation 
$(\partial_t - \Delta_f) u = \Sigma (t,x,u)$ verifying $0<u \le D$. Put $h = \log (u/D)$ and let $w = |\nabla h|^2/(1-h)^2$. 
Then $w$ satisfies the equation  
\begin{align}\label{eq21}
\square w = (\partial_t - \Delta_f) w =&-\frac{[\partial_t g + 2 {\mathscr Ric}_f]}{(1-h)^2} (\nabla h, \nabla h) 
- \frac{2h \langle \nabla h, \nabla w \rangle}{1-h} - 2 (1-h) w^2 \nonumber \\
&~- 2 \left| \frac{\nabla^2 h}{1-h} + \frac{\nabla h \otimes \nabla h}{(1-h)^2} \right|^2 
+ \frac{2 \langle \nabla h, \Sigma_x(t,x,De^h) \rangle}{D e^h (1-h)^2} \nonumber \\
&~+ 2 w \left[ \Sigma_u (t,x,De^h) + \frac{h \Sigma (t,x,De^h)}{De^h(1-h)} \right].
\end{align}
\end{lemma}

\begin{proof}
Referring to \eqref{h evolution eq} in Lemma \ref{lem20}, a straightforward calculation gives  
\begin{align}
\langle \nabla h, \nabla \partial_t h \rangle =&~\langle \nabla h, \nabla [\Delta_f h + |\nabla h|^2 + D^{-1} e^{-h} \Sigma(t,x,D e^h)] \rangle \\
=&~\langle \nabla h, \nabla \Delta_f h \rangle + \langle \nabla h, \nabla |\nabla h|^2 \rangle + |\nabla h|^2 \Sigma_u \nonumber \\
&~+ D^{-1} e^{-h} (- |\nabla h|^2 \Sigma + \langle \nabla h, \Sigma_x \rangle), \nonumber   
\end{align}
where we have abbreviated the arguments of $\Sigma$ and its partial derivatives $\Sigma_x, \Sigma_u$. Now by virtue of the identity 
$\partial_t |\nabla h|^2 = - [\partial_t g] (\nabla h, \nabla h) + 2 \langle \nabla h, \nabla \partial_t h \rangle$ 
({\it see} also \eqref{norm-grad-evolve-equation} in Lemma \ref{geometric-evolution-lemma-one}) it follows that 
\begin{align} \label{eq22a}
\partial_t |\nabla h|^2 =& \, - [\partial_t g] (\nabla h, \nabla h) + 2 \langle \nabla h, \nabla \Delta_f h \rangle +  2 \langle \nabla h, \nabla |\nabla h|^2 \rangle \nonumber \\
&+ \frac{2 \langle \nabla h, \Sigma_x \rangle}{De^h} + 2 |\nabla h|^2 \left( \Sigma_u - \frac{\Sigma}{De^h} \right).   
\end{align}
Moving next to the function $w$, by using $\partial_t w = [\partial_t |\nabla h|^2]/(1-h)^2 + [2|\nabla h|^2 \partial_t h]/(1-h)^3$, we have 
\begin{align}\label{eq22}
\partial_t w =& \, -\frac{[\partial_t g] (\nabla h, \nabla h)}{(1-h)^2} + \frac{2 \langle \nabla h, \nabla \Delta_f h \rangle}{(1-h)^2} 
+ \frac{2 \langle \nabla h, \nabla |\nabla h|^2 \rangle }{(1-h)^2} + \frac{2 \langle \nabla h, \Sigma_x \rangle}{D e^h (1-h)^2} \nonumber \\
&+ \frac{2 |\nabla h|^2}{(1-h)^2} \left( \Sigma_u - \frac{\Sigma}{De^h} \right)  + \frac{2 |\nabla h|^2  \Delta_f h}{(1-h)^3} + \frac{2 |\nabla h|^4}{(1-h)^3} 
+ \frac{2 |\nabla h|^2 \Sigma}{D e^h (1-h)^3}.
\end{align}
Likewise we have $\nabla w  = [\nabla |\nabla h|^2]/(1-h)^2 + [2 |\nabla h|^2 \nabla h]/(1-h)^3$ and so 
by recalling the relation $\Delta_f w  = \Delta w - \langle \nabla f, \nabla w \rangle$ it follows that 
\begin{align}\label{eq23}
\Delta_f w  =& \, \frac{\Delta_f |\nabla h|^2}{(1-h)^2} + \frac{4 \langle \nabla h, \nabla |\nabla h|^2 \rangle}{(1-h)^3} 
+ \frac{2 |\nabla h|^2 \Delta_f h}{(1-h)^3} + \frac{6 |\nabla h|^4}{(1-h)^4}. 
\end{align}
Putting (\ref{eq22})-(\ref{eq23}) together and taking into account the necessary cancellations gives  
\begin{align}
\Delta_f w 
=&~\partial_t w + \frac{[\partial_t g] (\nabla h, \nabla h)}{(1-h)^2}  
+ \frac{\Delta_f |\nabla h|^2}{(1-h)^2} 
- \frac{2 \langle \nabla h, \nabla \Delta_f h \rangle}{(1-h)^2} \nonumber \\
&~- \frac{2 \langle \nabla h, \nabla |\nabla h|^2 \rangle}{(1-h)^2}  
+ \frac{2 |\nabla h|^4}{(1-h)^3} 
+ \frac{4 \langle \nabla h, \nabla |\nabla h|^2 \rangle}{(1-h)^3} + \frac{6 |\nabla h|^4}{(1-h)^4} \nonumber \\
&~- \frac{2 \langle \nabla h, \Sigma_x \rangle}{D e^h (1-h)^2} 
- \frac{2 |\nabla h|^2}{(1-h)^2} \left( \Sigma_u - \frac{\Sigma}{De^h} \right) 
- \frac{2 |\nabla h|^2 \Sigma}{D e^h (1-h)^3}.
\end{align}
Now by making use of the weighted Bochner-Weitzenb\"ock formula \eqref{Bochner} this gives 
\begin{align}
\Delta_f w 
=&~\partial_t w + \frac{[\partial_t g] (\nabla h, \nabla h)}{(1-h)^2} 
+ \frac{2 |\nabla^2 h|^2}{(1-h)^2} + \frac{2 {\mathscr Ric}_f (\nabla h, \nabla h)}{(1-h)^2}  \nonumber \\
& - \frac{2 \langle \nabla h, \nabla |\nabla h|^2 \rangle}{(1-h)^2} 
+ \frac{4 \langle \nabla h, \nabla |\nabla h|^2 \rangle}{(1-h)^3} 
- \frac{2 |\nabla h|^4}{(1-h)^3}
+ \frac{6 |\nabla h|^4}{(1-h)^4} \nonumber \\
& - \frac{2 \langle \nabla h, \Sigma_x \rangle}{D e^h (1-h)^2} 
- \frac{2 |\nabla h|^2}{(1-h)^2} \left( \Sigma_u - \frac{\Sigma}{De^h} \right) 
- \frac{2 |\nabla h|^2 \Sigma}{D e^h (1-h)^3}, 
\end{align}
and therefore a rearrangement of terms and basic considerations leads to the formulation
\begin{align}
\Delta_f w 
=&~\partial_t w + \frac{[\partial_t g + 2 {\mathscr Ric}_f]}{(1-h)^2} (\nabla h, \nabla h) \nonumber \\
&+ 2 \left| \frac{\nabla^2 h}{1-h} + \frac{\nabla h \otimes \nabla h}{(1-h)^2} \right|^2 + \frac{2 |\nabla h|^4}{(1-h)^3} \nonumber \\
&- \frac{2 \langle \nabla h, \nabla |\nabla h|^2 \rangle}{(1-h)^2} 
- \frac{4 |\nabla h|^4}{(1-h)^3}+ \frac{2 \langle \nabla h, \nabla |\nabla h|^2 \rangle}{(1-h)^3} + \frac{4 |\nabla h|^4}{(1-h)^4} \nonumber \\
& - \frac{2 \langle \nabla h, \Sigma_x \rangle}{D e^h (1-h)^2} 
- \frac{2 |\nabla h|^2}{(1-h)^2} \left( \Sigma_u - \frac{\Sigma}{De^h} \right) 
- \frac{2 |\nabla h|^2 \Sigma}{D e^h (1-h)^3}.
\end{align}
Finally making note of the identity $(1-h)^3 \langle \nabla h, \nabla w \rangle=  2(1-h) \nabla^2 h(\nabla h, \nabla h) + 2 |\nabla h|^4$ results 
in the desired inequality. This therefore completes the proof. 
\end{proof}

\begin{lemma}\label{lem2112}
Under the assumption of Lemma $\ref{lem21}$, if the metric-potential pair $(g,f)$ evolves under the $\mathsf{k}$-super 
Perelman-Ricci flow inequality $\partial_t g + 2 {\mathscr Ric}_f(g) \ge -2 \mathsf{k} g$, 
then the function $w = |\nabla h|^2/(1-h)^2$ satisfies the parabolic differential inequality 
\begin{align}\label{eq2112}
\square w = (\partial_t - \Delta_f) w 
\le &-2(1-h)w^2 - \frac{2h \langle \nabla h, \nabla w \rangle}{1-h} \nonumber \\
&+ 2\mathsf{k}w + \frac{2 \langle \nabla h, \Sigma_x(t,x,De^h) \rangle}{D e^h (1-h)^2} \nonumber \\
&~+ 2 w \left[ \Sigma_u (t,x,De^h) + \frac{h \Sigma (t,x,De^h)}{De^h(1-h)} \right].
\end{align}
\end{lemma}

\begin{proof}
This is a straightforward consequence of the flow inequality \eqref{eq11c} and the identity \eqref{eq21} established 
in Lemma \ref{lem21} with the substitution $w = |\nabla h|^2/(1-h)^2$.
\end{proof}

\subsection{Localising in space-time and cut-offs}

In order to prove the local estimate in Theorem \ref{thm1}, we make use of the parabolic inequality 
in Lemma \ref{lem2112} in conjunction with a localisation argument. 
Towards this end fix $R, T>0$ and pick $\tau \in (0, T]$. Let $\varrho(x,t)$ denote the geodesic radial 
variable with respect to a fixed reference point $x_0$ at time $t$ and for $x \in M$ and $0 \le t \le T$ set  
\begin{equation} \label{cut-off def}
\phi(x,t) = \bar{\phi}(\varrho(x,t), t).
\end{equation}
Here $\bar\phi$ is a suitable function of the real non-negative variables $\varrho$ and $t$. 
The resulting space-time function $\phi$ will then serve as a smooth cut-off function supported in the compact set 
$Q_{R,T} \subset M \times [0, T]$. The existence of $\bar{\phi}=\bar \phi(\varrho, t)$ as 
used in \eqref{cut-off def} and its properties is granted by the following result (see \cite{[Ba1], Bri, SZ, [Wu15]}).

\begin{lemma} \label{phi lemma} There exists a smooth function $\bar{\phi}:[0,\infty)\times [0,T] \to \mathbb{R}$ such that:
\begin{enumerate}[label=$(\roman*)$]
\item ${\rm supp} \, \bar{\phi}(\varrho,t) \subset [0,R] \times [0,T]$  and $0\leq \bar{\phi}(\varrho,t) \leq 1$ in $[0,R] \times [0,T]$.
\item $\bar{\phi}=1$ in $[0,R/2] \times [\tau, T]$ and $\partial \bar{\phi}/\partial \varrho =0$ in $[0,R/2] \times [0, T]$, respectively.
\item $\bar{\phi}(\varrho,0)=0$ for all $\varrho \in [0,\infty)$ and there exists $c>0$ such that the bound
\begin{equation}
\left| \frac{\partial \bar{\phi}}{\partial t} \right| \le c \frac{\bar{\phi}^{1/2}}{\tau}, 
\end{equation} 
holds on $[0,\infty)\times[0,T]$.
\item For every $0<a<1$ there exists $c_a>0$ such that the bounds 
\begin{equation}
-c_a \frac{\bar{\phi}^a}{R} \le \frac{\partial \bar{\phi}}{\partial \varrho} \le 0, \qquad 
\left|\frac{\partial^2 \bar{\phi}}{\partial \varrho^2}\right| \le c_a \frac{\bar{\phi}^a}{R^2},  
\end{equation}
hold on $[0, \infty)\times [0, T]$.
\end{enumerate}
\end{lemma}

\subsection{Proof of the local estimate in Theorem \ref{thm1}}

This will be carried out in two stages. First we establish the estimate in the case $q \equiv 0$ and then we pass to the general case. 
So for now suppose $q \equiv 0$. Fix $\tau \in (0, T]$ and $\phi(x,t) = \bar{\phi}(\varrho(x,t), t)$ with $\bar \phi$ as in 
Lemma \ref{phi lemma}. We show the respective estimate to hold at $(x,\tau)$ with $d(x, x_0, \tau) \le R/2$. The 
arbitrariness of $\tau>0$ will then give the estimate for all $(x, t)$ in $Q_{R/2, T}$ satisfying $t>0$. Now starting with 
the localised function $\phi w$ it is clear that  
\begin{equation}
\square (\phi w) = \phi \square w 
- 2 [\langle\nabla \phi ,\nabla (\phi w) \rangle- |\nabla \phi|^2 w]/\phi  + w \square \phi.
\end{equation}  
Substituting from \eqref{eq2112} in Lemma \ref{lem2112} (whilst recalling the relation $\mathsf{k} = (n-1) k_1 + k_2$) then gives
\begin{align} \label{eq25}
\square (\phi w) = (\partial_t - \Delta_f) (\phi w) 
\le &-\left\langle \frac{2h\nabla h}{1 -h} + \frac{2\nabla \phi}{\phi}, \nabla (\phi w) \right\rangle \nonumber \\
&+ w \left\langle \frac{2h \nabla h}{1 -h} + \frac{2 \nabla \phi}{\phi}, \nabla \phi \right\rangle -2(1-h) \phi w^2 \nonumber \\
&~+ w (\partial_t - \Delta_f + 2\mathsf{k}) \phi + \frac{2 \phi \langle \nabla h, \Sigma_x(t,x,De^h) \rangle}{D e^h (1-h)^2} \nonumber \\
&~+2\phi w \left[ \Sigma_u (t,x,De^h) + \frac{h \Sigma (t,x,De^h)}{De^h(1-h)} \right]. 
\end{align}

Suppose now that the localised function $\phi w$ is maximised at the point $(x_1, t_1)$ in the compact set 
$\{(x,t) : d(x,x_0, t) \le R \mbox{ and } 0 \le t \le \tau\} \subset M \times [0, T]$. By virtue of Calabi's standard argument 
({\it cf.} \cite{Calabi} or \cite{SchYau} p.~21) we can assume that $x_1$ is not in 
the cut locus of $x_0$ and so $\phi$ is smooth at $(x_1,t_1)$ for the application of the maximum principle. Additionally, we can assume that $(\phi w)(x_1, t_1) >0$ 
as otherwise the conclusion is true with $w(x, \tau) \le 0$ for all $d(x, x_0, \tau) \le R/2$. In particular 
$t_1>0$ and at the point $(x_1,t_1)$  we have $\partial_t (\phi w) \ge 0$, $\nabla(\phi w) =0$ and $\Delta_f(\phi w) \le 0$. From \eqref{eq25} it thus follows that  
\begin{align} 
2(1-h)\phi w^2 
\le &~w \left\langle \frac{2h \nabla h}{1 -h} + \frac{2 \nabla \phi}{\phi}, \nabla \phi \right\rangle
- w \Delta_f \phi + w \partial_t \phi  + 2\mathsf{k}w\phi \nonumber \\
&~+ \frac{2 \phi \langle \nabla h, \Sigma_x(t,x,De^h) \rangle}{D e^h (1-h)^2} 
+ 2 \phi w \left[ \Sigma_u (t,x,De^h) + \frac{h \Sigma (t,x,De^h)}{De^h(1-h)} \right] \nonumber 
\end{align}
at $(x_1,t_1)$ or dividing through by $2(1-h) \ge 2$ that 
\begin{align}\label{eq27}
\phi w^2 
\le &~w \left\langle \frac{h \nabla h}{1 -h} + \frac{\nabla \phi}{\phi}, \frac{\nabla \phi}{1-h} \right\rangle
+ \frac{w[-\Delta_f \phi + \partial_t \phi + 2\mathsf{k}\phi]}{2(1-h)} \nonumber \\
&~+ \frac{\phi \langle \nabla h, \Sigma_x(t,x,De^h) \rangle}{D e^h (1-h)^3} 
+ \phi w \left[ \frac{\Sigma_u (t,x,De^h)}{1-h} + \frac{h \Sigma (t,x,De^h)}{De^h(1-h)^2} \right].
\end{align}
The goal is now to use \eqref{eq27} to establish the required estimate at $(x, \tau)$. To this end we 
consider two cases. Firstly, we consider $d(x_1, x_0, t_1) \le 1$ and next  
$d(x_1, x_0, t_1) \ge 1$. \\
{\bf Case 1.} Since here $\phi$ is a constant 
function in the space direction [for all $x$ satisfying $d(x, x_0,t) \le R/2$, where $t \in [0, T]$ and $R \geq 2$ by property 
$(ii)$] all the terms involving space derivatives of $\phi$ at $(x_1, t_1)$ vanish (in particular $\nabla \phi=0$, 
$\Delta_f \phi =0$ and $\phi_t = \bar \phi_t$). So as a result it follows from (\ref{eq27}) that at 
the point $(x_1,t_1)$, we have the bound
\begin{align*}
\phi w^2 
\le \frac{w[|\phi_t|+2\mathsf{k}\phi]}{2(1-h)} + \frac{\phi |\langle \nabla h, \Sigma_x(t,x,De^h) \rangle|}{D e^h (1-h)^3} 
+ \phi w \left[ \frac{\Sigma_u (t,x,De^h)}{1-h} + \frac{h \Sigma (t,x,De^h)}{De^h(1-h)^2} \right]_+
\end{align*}

Now upon writing  
$w[|\phi_t|+2\mathsf{k}\phi] = \sqrt \phi w [|\phi_t|/\sqrt \phi+2\mathsf{k}\sqrt \phi] \le \sqrt \phi w [c/\tau+2\mathsf{k}\sqrt \phi]$, 
hence giving $w[|\phi_t|+2\mathsf{k}\phi]/[2(1-h)] \le (1/4) \phi w^2 + C [1/\tau^2 + \mathsf{k}^2]/(1-h)^2$ (by an application of 
the Cauchy-Schwarz inequality and $(iii)$ in Lemma \ref{phi lemma}) and using similar bounds on the other two terms 
on the right-hand side by utilising Young's inequality, we obtain after rearranging terms, 
\begin{align*}
\phi w^2 \le C \left\{ \frac{1/\tau^2 + \mathsf{k}^2}{(1-h)^2} 
+ \left[ \frac{\Sigma_u (t,x,De^h)}{1-h} + \frac{h \Sigma (t,x,De^h)}{De^h(1-h)^2} \right]^2_+ 
+ \left[ \frac{|\Sigma_x(t,x,De^h)|}{D e^h (1-h)^2} \right]^{4/3} \right\} . \nonumber
\end{align*}

As $\phi \equiv 1$, when $d(x,x_0,t) \le R/2$ and $\tau \le t \le T$ (hence in particular for $t = \tau$) 
by $(ii)$ in Lemma \ref{phi lemma}, we have $w(x,\tau) =(\phi w)(x, \tau) \leq (\phi w)(x_1,t_1) \le (\sqrt \phi w) (x_1,t_1)$. Thus recalling 
the definitions of $\mathsf{k}$, $k$, $w=|\nabla h|^2/(1-h)^2$, $h=\log(u/D)$, noting $1/(1-h)\le1$, and adjusting the 
constant if necessary, we arrive at the bound at $(x, \tau)$ 
\begin{align*}
\frac{\left| \nabla h \right|}{1- h} \le C \left\{ \frac{1}{\sqrt \tau} + \sqrt k 
+ \sup_{Q_{R,T}} \left[ \frac{De^h (1-h) \Sigma_u + h \Sigma}{De^h(1-h)^2} \right]_+^{1/2} 
+ \sup_{Q_{R,T}} \left[ \frac{|\Sigma_x|}{De^h (1-h)^2} \right]^{1/3} \right\}.   
\end{align*}
This together with the arbitrariness of $\tau>0$ is now immediately seen to be a special case of the estimate \eqref{eq13}.

\qquad \\
{\bf Case 2.} Upon referring to the right-hand side of \eqref{eq27}, and noting the properties of $\bar \phi$ as listed in Lemma \ref{phi lemma} 
we proceed onto bounding the full expression on the right-hand side on \eqref{eq27} in the case $d(x_1, x_0, t_1) \ge 1$. Towards this end 
dealing with the first term first, we have 
\begin{align} \label{eq28}
\left\langle w \left[ \frac{h \nabla h}{1 -h} + \frac{\nabla \phi}{\phi} \right], \frac{\nabla \phi}{1-h} \right\rangle
&\le w \left[ \frac{h |\nabla h|}{1 -h} + \frac{|\nabla \phi|}{\phi} \right] \frac{|\nabla \phi|}{1-h} \nonumber \\
&\le w \left[ h \sqrt w + \frac{|\nabla \phi|}{\phi} \right] \frac{|\nabla \phi|}{1-h} \nonumber \\
&\le w \sqrt \phi \left [ \phi^{1/4} \frac{\sqrt w |h|}{1-h} \frac{|\nabla \phi|}{\phi^{3/4}} + \frac{|\nabla \phi|^2}{\phi^{3/2}} \right] \nonumber \\
&\le \frac{\phi w^2}{4} + \frac{C}{R^4}. 
\end{align}
In much the same way regarding the terms involving $\Sigma$ we have firstly 
\begin{align}
\frac{\phi \langle \nabla h, \Sigma_x \rangle}{D e^h (1-h)^3} &\le \frac{\phi |\nabla h| |\Sigma_x|}{D e^h (1-h)^3} 
= \frac{\phi \sqrt w |\Sigma_x|}{D e^h (1-h)^2} \nonumber \\
&\le \frac{\phi w^2}{8} + C \left(\frac{|\Sigma_x|}{De^h (1-h)^2} \right)^{4/3}
= \frac{\phi w^2}{8} + C \mathsf{P}_\Sigma^{4/3}, 
\end{align}
with $\mathsf{P}_\Sigma=|\Sigma_x|/[(1-h)^2 De^h]$ and likewise for the subsequent terms, upon noting $-1 \le h/(1-h) \le 0$, $h \le 0$ and $0 \le \phi \le 1$, we have
\begin{align}
\phi w \left( \frac{\Sigma_u}{1-h} + \frac{h \Sigma}{De^h(1-h)^2} \right) &= \phi w \left( \frac{De^h (1-h) \Sigma_u + h \Sigma}{De^h(1-h)^2} \right) \nonumber \\
&\le \frac{\phi w^2}{8} + C \left[ \frac{De^h (1-h) \Sigma_u + h \Sigma}{(1-h)^2 De^h} \right]^2_+ \nonumber \\
&= \frac{\phi w^2}{8} + C \mathsf{R}^2_\Sigma, 
\end{align}
where in the last equation we have set $\mathsf{R}_\Sigma = \{ [De^h (1-h) \Sigma_u + h \Sigma]/[(1-h)^2 De^h] \}_+$.

Now for the $\Delta_f \phi$ term we use the Wei-Wylie weighted Laplacian comparison theorem taking advantage of the fact that it only depends 
on the lower bound on ${\mathscr Ric}_f(g)$ (\cite{[WeW09]}). Indeed recalling $\varrho(x, t)=d(x, x_0, t)$, $1 \le d(x_1,x_0, t_1) \le R$, it follows from 
Theorem 3.1 in \cite{[WeW09]} and ${\mathscr Ric}_f (g) \ge -(n-1) k_1 g$ with $k_1 \ge 0$, that 
\begin{equation}
\Delta_f \varrho \le \gamma_{\Delta_f}+(n-1)(R-1) k_1,
\end{equation}
whenever $1 \le \varrho \le R$, $0 \le t \le T$ [therefore in particular at the space-time point $(x_1, t_1)$]. 
Here as indicated earlier we have set 
\begin{equation}
\gamma_{\Delta_f} = \max_{(x,t)}  \{ \Delta_f \varrho (x,t) : d(x,x_0,t)=1, \, 0 \le t \le T \}.
\end{equation} 
Thus proceeding on to bounding $-\Delta_f \phi$, upon referring to \eqref{cut-off def} and using $(ii)$ [$\bar \phi_r =0$ 
when $0 \le r \le R/2$], $(iv)$ [$\bar \phi_\varrho \le 0$ when $0 \le \varrho < \infty$ together with the bounds] 
we have $-\Delta_f \phi = -(\bar \phi_{\varrho \varrho} |\nabla \varrho|^2 + \bar \phi_\varrho \Delta_f \varrho)$ and so 
\begin{align} \label{bound Delta f alpha eq}
- \Delta_f \phi &\le \left( \frac{|\bar \phi_{\varrho \varrho}|}{\sqrt{\bar \phi}} + \frac{|\bar \phi_\varrho|}{\sqrt{\bar \phi}} 
\left( [\gamma_{\Delta_f}]_+ +(n-1) (R-1) k_1 \right) \right) \sqrt{\bar \phi} \nonumber \\
&\le C \left(  \frac{1}{R^2} + \frac{[\gamma_{\Delta_f}]_+}{R} + k_1 \right) \sqrt{\phi}. 
\end{align}

Next to bound the term $\partial_t \phi$ pick $x$ such that $d(x, x_0; t) \le R$ and let $\gamma: [0,1] \to M$ 
be a minimal geodesic connecting $x_0$ and $x$ at the fixed time $t$ where we write $\gamma=\gamma(s)$ 
with $\gamma(0)=x_0$ and $\gamma(1)=x$. Then, recalling 
the bound $\partial_t g \ge -2 k_2 g$ in $\mathcal{B}_{R, T}$, we have 
\begin{align}
\partial_t \varrho (x,t) &= \partial_t \int_0^1 |\gamma'(s)|_{g(t)} \, ds = \int_0^1 \partial_t [|\gamma'(s)|_{g(t)}] \, ds \nonumber \\
&= \int_0^1 \frac{[\partial_t g] (\gamma', \gamma')}{2 |\gamma'|_{g(t)}} \, ds \nonumber \\
&\ge \int_0^1 - k_2 |\gamma'|_{g(t)} \, ds 
= - k_2 \varrho \ge - k_2 R.
\end{align}
Hence a straightforward differentiation followed by an application of the properties of $\bar \phi$ in Lemma \ref{phi lemma} gives 
$\partial_t \phi = \bar{\phi}_t + \bar \phi_\varrho \partial_t \varrho \le |\bar \phi_t| - k_2 R \bar \phi_\varrho 
= |\bar \phi_t| + k_2 R |\bar \phi_\varrho| \le C [ 1/\tau + k_2] \sqrt{\phi}$. 
Therefore, by summarising, the above estimates on the derivatives of $\phi$ along with the term $2\mathsf{k} \phi w$ in \eqref{eq27} give, 
after an application of Young inequality, the bounds
\begin{align} 
& \frac{w (-\Delta_f) \phi}{2(1-h)} \le \frac{\phi w^2}{8} + C \left( \frac{1}{R^4} + \frac{[\gamma_{\Delta_f}]_+^2}{R^2} + k_1^2 \right), \label{eq28b} \\ 
& \frac{w (2\mathsf{k} + \partial_t) \phi}{2(1-h)} \le \frac{\phi w^2}{8} + C \left( \frac{1}{\tau^2} + \mathsf{k}^2 + k_2^2 \right). \label{eq211}
\end{align}

Now referring to \eqref{eq27}, noting the inequality $1-h \ge 1$ and making use of the bounds obtained in \eqref{eq28}--\eqref{eq211}, 
it follows after reverting to $u=De^h$, writing $k=(k_1^2+k_2^2)^{1/2}$ and noting 
$\mathsf{k}=(n-1)k_1+k_2$, that the following upper bound holds for $\phi w^2$ at $(x_1, t_1)$,  
\begin{align}
\phi w^2 \le&  \, C \left\{ \frac{1+R^2[\gamma_{\Delta_f}]^2_+}{R^4} + \frac{1}{\tau^2} + k^2 
+ \sup_{Q_{R,T}} \left( \mathsf{R}^{2}_\Sigma(u) + \mathsf{P}^{4/3}_\Sigma(u) \right) \right\}.
\end{align}

Recalling the maximality of $\phi w$ at $(x_1, t_1)$ along with $\phi \equiv 1$ when $d(x, x_0; t) \le R/2$ and $\tau \le t \le T$, 
it follows that $w^2 (x, \tau) = (\phi^2 w^2)(x, \tau) \le (\phi^2 w^2)(x_1,t_1) \leq (\phi w^2)(x_1,t_1)$ 
when $d(x, x_0; \tau) \le R/2$. Hence upon noting $w = |\nabla h|^2/(1-h)^2$, the above gives
\begin{align}\label{eq214}
\frac{\left| \nabla \log u \right|}{1- \log(u/D)} 
\le& C \left\{ \frac{1}{R} +\sqrt{\frac{[\gamma_{\Delta_f}]_+}{R}} + \frac{1}{\sqrt{\tau}} + \sqrt{k} 
+ \sup_{Q_{R,T}} \left( \mathsf{R}^{1/2}_\Sigma(u) + \mathsf{P}^{1/3}_\Sigma(u) \right)
 \right\}.
\end{align}

Thus in either case the estimate is true at $(x, \tau)$ and so the desired conclusion when $q \equiv 0$ and hence 
$\mathsf{N} \equiv 0$ follows from the arbitrariness of $\tau \in (0, T]$.

The passage to the general case (non-zero $q$) now follows by replacing $\Sigma$ with $\Sigma+qu$. Indeed 
by the sub-additivity of $[\cdot]_+$ we have $\mathsf{R}_{\Sigma + qu} \le \mathsf{R}_{\Sigma} + \mathsf{R}_{qu}$ and  
$\mathsf{P}_{\Sigma + qu} \le \mathsf{P}_{\Sigma} + \mathsf{P}_{qu}$. Moreover, referring to \eqref{eq13-R} it is easily seen that,  
\begin{equation}
\mathsf{R}_{qu} = \left[ \frac{q}{1-h} + \frac{hqu}{u(1-h)^2} \right]_+ = \left[ \frac{qu(1-h) + hqu}{u(1-h)^2} \right]_+ = \frac{q_+}{(1-h)^2} \le q_+,
\end{equation}
and 
\begin{equation}
\mathsf{P}_{qu} = \frac{|u \nabla q|}{u(1-h)^2} = \frac{|\nabla q|}{(1-h)^2} \le |\nabla q|.
\end{equation}
Thus invoking the conclusion of the theorem from the first part (when $q \equiv 0$) and using the above calculation, the desired 
estimate \eqref{eq13} for general $q$ follows by noting,  
\begin{align}
\mathsf{R}^{1/2}_{\Sigma+qu} + \mathsf{P}^{1/3}_{\Sigma+qu}
&\le [\mathsf{R}_{\Sigma} + \mathsf{R}_{qu}]^{1/2} + [\mathsf{P}_{\Sigma} + \mathsf{P}_{qu}]^{1/3} \nonumber \\
&\le \mathsf{R}^{1/2}_{\Sigma} + \mathsf{R}^{1/2}_{qu}  +  \mathsf{P}^{1/3}_{\Sigma} + \mathsf{P}^{1/3}_{qu} \nonumber \\
&\le \mathsf{R}^{1/2}_{\Sigma} + q_+^{1/2} + \mathsf{P}^{1/3}_{\Sigma} + |\nabla q|^{1/3}
= \mathsf{N}_q + \mathsf{R}^{1/2}_{\Sigma} + \mathsf{P}^{1/3}_{\Sigma}. 
\end{align}
The proof is thus complete. \hfill $\square$

\subsection{Proof of the elliptic Harnack inequality in Theorem \ref{cor Harnack}} 

We now come to the proof of the Harnack inequality as formulated in Theorem \ref{cor Harnack}. In its local form this is a consequence 
of the local estimate in Theorem \ref{thm1} and in its global form this is a consequence of the global estimate in Theorem 
\ref{thm1-global}. Below we give the proof in the local case as the other is similar.

Towards this end fix $t \in (0, T)$ and pick $(x_1,t)$ and $(x_2,t)$ in $Q_{R/2, T}$ as in the theorem. Let $\zeta=\zeta(s)$ with $0 \le s \le 1$ 
be a shortest curve with respect to $g=g(t)$ joining $x_1$, $x_2$ [that is, $\zeta(0)=x_1$, $\zeta(1)=x_2$] with $(\gamma(t), t)$ lying 
entirely in $Q_{R/2, T}$. Put $d=d(x_1, x_2, t)$ and 
assume that $\mathsf{R}_\Sigma, \mathsf{P}_\Sigma < \infty$ (otherwise $\gamma=0$ and the inequality is trivially true). 
Now by utilising the estimate \eqref{eq13} in Theorem \ref{thm1} we can write 
\begin{align}
\log \frac{1-h(x_2, t)}{1-h(x_1, t)} &= \int_0^1 \frac{d}{ds} \log [1-h(\zeta(s), t)] \, ds 
= \int_0^1 - \frac{\langle \nabla h (\zeta(s), t), \zeta'(s) \rangle}{1-h(\zeta(s), t)} \, ds \nonumber \\
&\le \int_0^1 \frac{|\nabla h||\zeta'|}{1-h} \, ds \le \left[ \sup_{Q_{R/2,T}} \frac{|\nabla h|}{1-h} \right] \int_0^1 |\zeta'| \, ds \nonumber \\
&\le C \left\{ \frac{1}{R} + \frac{1}{\sqrt{t}} + \sqrt{k} 
+ \sup_{Q_{R,T}} \left[ \mathsf{N}_q + \mathsf{R}^{1/2}_\Sigma 
+ \mathsf{P}^{1/3}_\Sigma \right] 
+ \sqrt{\frac{[\gamma_{\Delta_f}]_+}{R}} 
\right\}. 
\end{align}
Here we have used $|\nabla h|/(1-h) = |\nabla \log u|/[1-\log(u/D)]$. Therefore a direct set of calculations give, 
\begin{align}
\frac{1-h(x_2, t)}{1-h(x_1, t)}
&= \frac{1 - \log [u(x_2, t)/D]}{1 - \log [u(x_1, t)/D]}
=\frac{\log [eD/u(x_2, t)]}{\log [eD/u(x_1, t)]} \\
&\le {\rm exp} \left[ C d \left( \frac{1}{R} 
+ \sqrt k + \frac{1}{\sqrt{t}}  
+ \sup_{Q_{R,T}} \left[ \mathsf{N}_q + \mathsf{R}^{1/2}_\Sigma 
+ \mathsf{P}^{1/3}_\Sigma \right] 
+ \sqrt{\frac{[\gamma_{\Delta_f}]_+}{R}}
\right) \right] \nonumber 
\end{align}
and so the assertion follows. As indicated earlier the proof of the global version of the inequality is similar and 
is hence abbreviated. \hfill $\square$

\section{Proof of the elliptic Hamilton estimate in Theorem \ref{thm18}}\label{sec5}

\subsection{Some intermediate parabolic lemmas (${\bf II}$)}

Before moving onto presenting the proof of Theorem \ref{thm18} we gather together some useful components and 
tools that are needed later. In the spirit of what was done earlier in Section \ref{sec3} we present here a set of parabolic 
identities for suitable quantities built out the solution.

\begin{lemma} \label{lemma h one}
Suppose $u$ is a positive solution to $(\partial_t - \Delta_f) u = \Sigma (t,x,u)$ and put $h=u^\beta$ where $\beta \in (0, 1)$. 
Then $h$ satisfies the equation 
\begin{align} \label{eq4.1}
\square h = (\partial_t - \Delta_f) h = (1-\beta)|\nabla h|^2/(\beta h) + \beta h^{1-1/\beta} \Sigma(t,x, h^{1/\beta}).
\end{align}
\end{lemma}

\begin{proof}
A basic calculation gives $\Delta_f h = \beta (\beta-1) u^{\beta-2} |\nabla u|^2 + \beta u^{\beta -1} \Delta_f u$ and 
$\partial_t h = \beta u^{\beta-1} \partial_t u$. Therefore \eqref{eq4.1} follows at once by substitution in \eqref{eq11}. 
\end{proof}

\begin{lemma} \label{lemma h two}
Under the assumptions of Lemma $\ref{lemma h one}$ on $u$, $h$, the function $|\nabla h|^2$ satisfies the equation 
\begin{align}
\square |\nabla h|^2 =  (\partial_t - \Delta_f) |\nabla h|^2   
=&-[\partial_t g (\nabla h, \nabla h)+ 2 {\mathscr Ric}_f(\nabla h, \nabla h)] \nonumber \\
&-2 |\nabla^2 h|^2 + \frac{2(\beta-1)}{\beta h^2} \left[ |\nabla h|^4 - h \langle \nabla h, \nabla |\nabla h|^2 \rangle \right] \nonumber \\
&+2 \beta \langle \nabla h, \nabla [h^{1-1/\beta} \Sigma(t,x,h^{1/\beta})] \rangle. 
\end{align}
\end{lemma}

\begin{proof}
By making use of the weighted Bochner-Weitzenb\"ock formula \eqref{Bochner} and the evolution of $h$ described by \eqref{eq4.1} we have 
\begin{align*}
(\partial_t - \Delta_f) |\nabla h|^2  
=&~\partial_t |\nabla h|^2 - 2 |\nabla^2 h|^2 - 2 \langle \nabla h, \nabla \Delta_f h \rangle - 2 {\mathscr Ric}_f(\nabla h, \nabla h) \nonumber \\
=&~- g_t (\nabla h, \nabla h) + 2 \langle \nabla h, \nabla (\partial_t - \Delta_f) h \rangle - 2 |\nabla^2 h|^2 - 2 {\mathscr Ric}_f(\nabla h, \nabla h)  \nonumber \\
=&~-[g_t (\nabla h, \nabla h) + 2 {\mathscr Ric}_f(\nabla h, \nabla h)] - 2 |\nabla^2 h|^2 \nonumber \\
&~+ \langle 2(1-\beta) \nabla h, \nabla [|\nabla h|^2/(\beta h)] \rangle 
+ 2 \beta \langle \nabla h, \nabla [h^{1-1/\beta} \Sigma (t,x,h^{1/\beta})]  \rangle. 
\end{align*}
Expanding the gradient and the inner product then gives the desired conclusion.
\end{proof}

With the conclusions of the above two lemmas at hand let us now calculate $\square (h |\nabla h|^2)$. 
To this end first using the weighted Bochner-Weitzenb\"ock formula we can write 
\begin{align*} 
\Delta_f (h |\nabla h|^2) =&~h \Delta_f |\nabla h|^2 + |\nabla h|^2\Delta_f h + 2 \langle \nabla h, \nabla |\nabla h|^2 \rangle \nonumber \\
=&~2h |\nabla^2 h|^2 + 2 h \langle \nabla h, \nabla \Delta_f h \rangle + 2h {\mathscr Ric}_f(\nabla h, \nabla h) \nonumber \\
&+ |\nabla h|^2 \Delta_f h + 2 \langle \nabla h, \nabla |\nabla h|^2 \rangle, 
\end{align*}
and next we have 
\begin{align*}
\partial_t (h |\nabla h|^2) &= |\nabla h|^2 \partial_t h + h \partial_t |\nabla h|^2 
= |\nabla h|^2 \partial_t h - h [\partial_t g] (\nabla h, \nabla h) + 2h \langle \nabla h, \nabla \partial_t h \rangle.
\end{align*} 
Thus putting the two together gives 
\begin{align} \label{heatfW thm18}
\square (h |\nabla h|^2) =&-2h |\nabla^2 h|^2 + 2 h \langle \nabla h, \nabla \square h \rangle - 2h {\mathscr Ric}_f(\nabla h, \nabla h) \nonumber \\
&~+ |\nabla h|^2 \square h - 2 \langle \nabla h, \nabla |\nabla h|^2 \rangle - h [\partial_t g](\nabla h, \nabla h). 
\end{align}

\subsection{Proof of the local estimate in Theorem \ref{thm18}}
We consider first the case $q \equiv 0$. From \eqref{eq11c} and $h \ge 0$ we have  
$h [\partial_t g](\nabla h, \nabla h) + 2h {\mathscr Ric}_f(\nabla h, \nabla h) \ge - 2\mathsf{k}h |\nabla h|^2$ and so substituting into 
\eqref{heatfW thm18} and making use of Lemma \ref{lemma h one} results in 
\begin{align} \label{Pf thm18 eq1}
(\partial_t - \Delta_f) [h |\nabla h|^2] \le &-2h |\nabla^2 h|^2 - 2(\beta-1)
[h \langle \nabla h, \nabla |\nabla h|^2 \rangle - |\nabla h|^4]/(\beta h) + 2\mathsf{k}h |\nabla h|^2 \nonumber \\
&~+2 \beta h \langle \nabla h, \nabla [h^{1-1/\beta} \Sigma(t,x, h^{1/\beta})] \rangle - (\beta-1) |\nabla h|^4/ (\beta h) \nonumber \\
&~+\beta h^{1-1/\beta} \Sigma(t,x, h^{1/\beta}) |\nabla h|^2 - 2 \langle \nabla h, \nabla |\nabla h|^2 \rangle. 
\end{align}
Using $2|\sqrt{h} \nabla^2 h + [\nabla h \otimes \nabla h]/\sqrt{h}|^2 
= 2 h |\nabla^2 h|^2 + 2 \langle \nabla h, \nabla |\nabla h|^2 \rangle + 2|\nabla h|^4/h \ge 0$ 
we can rewrite \eqref{Pf thm18 eq1} as
\begin{align} \label{Pf thm18 eq2}
(\partial_t - \Delta_f) [h |\nabla h|^2] 
\le&-[2(\beta-1)/\beta] \langle \nabla h, \nabla |\nabla h|^2 \rangle  
+ [(3\beta-1)/(\beta h)] |\nabla h|^4 \nonumber \\
&+2\mathsf{k} h |\nabla h|^2 + 2 \beta h \langle \nabla h, \nabla [h^{1-1/\beta} \Sigma(t,x, h^{1/\beta})] \rangle \nonumber \\
&+\beta h^{1-1/\beta} \Sigma(t,x, h^{1/\beta}) |\nabla h|^2,  
\end{align}
and upon making note of $\langle \nabla h, \nabla (h |\nabla h|^2) \rangle = |\nabla h|^4 + h \langle \nabla h, \nabla |\nabla h|^2 \rangle$ 
and substituting in \eqref{Pf thm18 eq2} we can write  
\begin{align}
(\partial_t - \Delta_f) [h |\nabla h|^2] 
\le &-[2(\beta-1)/(\beta h)] \langle \nabla h, \nabla (h|\nabla h|^2) \rangle \nonumber \\
&-[(3-5\beta)/(\beta h^3)] (h |\nabla h|^2)^2 + 2\mathsf{k} (h |\nabla h|^2)  \nonumber \\
&+2\beta h \langle \nabla h, \nabla [h^{1-1/\beta} \Sigma(t,x, h^{1/\beta})] \rangle \nonumber \\
&+\beta h^{1-1/\beta} \Sigma(t,x, h^{1/\beta}) |\nabla h|^2. 
\end{align}

Let $\mathsf{Z}_\Sigma$ denote the sum of the last two terms on the right-hand side of \eqref{Pf thm18 eq2} 
and let us abbreviate the arguments of $\Sigma, \Sigma_x, \Sigma_u$ for convenience. A basic calculation gives  
\begin{align} \label{Pf thm18 eq3}
\langle \nabla h, \nabla [h^{1-1/\beta} \Sigma] \rangle 
= \frac{(\beta-1)}{\beta h^{1/\beta}} |\nabla h|^2 \Sigma + h^{1-1/\beta} \langle \nabla h, \nabla \Sigma \rangle, 
\end{align}
where the last term on the right in \eqref{Pf thm18 eq3} can in turn be calculated as 
\begin{align} \label{Pf thm18 eq4}
\langle \nabla h, \nabla \Sigma \rangle &= \langle \nabla h, \Sigma_x + \nabla (h^{1/\beta}) \Sigma_u \rangle \nonumber \\
&= \langle \nabla h, \Sigma_x \rangle + \beta^{-1} h^{1/\beta -1} |\nabla h|^2 \Sigma_u \nonumber \\
&= \langle \nabla h, \Sigma_x \rangle + \beta^{-1} h^{1/\beta -2} (h |\nabla h|^2) \Sigma_u. 
\end{align}
Substituting \eqref{Pf thm18 eq3}-\eqref{Pf thm18 eq4} back into the expression for $\mathsf{Z}_\Sigma$ leads to
\begin{align}
\mathsf{Z}_\Sigma &= 2 \beta h \langle \nabla h, \nabla [h^{1-1/\beta} \Sigma] \rangle 
+ \beta h^{1-1/\beta} |\nabla h|^2 \Sigma \nonumber \\
& = (3\beta-2) h^{1-1/\beta} |\nabla h|^2 \Sigma + 2\beta h^{2-1/\beta} \langle \nabla h, \nabla \Sigma \rangle \nonumber \\
&= \frac{(3\beta-2) \Sigma + 2u \Sigma_u}{u} (h |\nabla h|^2) + \frac{2\beta h^2 \langle \nabla h, \Sigma_x \rangle}{u}.
\end{align}

We now take a space-time cut-off function $\phi$ as in \eqref{cut-off def} and consider the localised function $\phi h|\nabla h|^2$. 
Then as in the proof of Theorem \ref{thm1} we can write 
\begin{align}  \label{Pf thm18 eq5}
(\partial_t-\Delta_f) [\phi h|\nabla h|^2]
\le &-\left\langle 2 \left[ \frac{(\beta-1) \nabla h}{\beta h} + \frac{\nabla \phi}{\phi} \right], \nabla(\phi h|\nabla h|^2) \right\rangle \nonumber \\
&~+2 (h|\nabla h|^2) \left\langle \frac{(\beta-1) \nabla h}{\beta h} + \frac{\nabla \phi}{\phi}, \nabla \phi \right\rangle\\
&~- \frac{3-5\beta}{\beta h^3} \phi (h|\nabla h|^2)^2 + h|\nabla h|^2 (\partial_t - \Delta_f + 2\mathsf{k}) \phi 
+ \mathsf{Z}_\Sigma \phi. \nonumber 
\end{align}

For fixed $\tau \in (0,T]$ let $(x_1, t_1)$ be a maximum point for the localised function $\phi h|\nabla h|^2$ 
in the compact set $\{d(x,x_0,t) \le R, 0 \le t \le \tau\} \subset M \times [0, T]$. 
Without loss of generality we can take $t_1>0$ and for the sake of establishing the estimate 
at $(x, \tau)$ in $Q_{R/2, T}$ it suffices to confine to the case $d(x_1, x_0, t_1) \ge 1$. Now at 
$(x_1, t_1)$ we have the relations 
$\partial_t (\phi h|\nabla h|^2) \ge 0$, $\nabla(\phi h|\nabla h|^2) =0$ and $\Delta_f(\phi h|\nabla h|^2) \le 0$. 
Therefore applying these to \eqref{Pf thm18 eq5} and rearranging the inequality result in 
\begin{align} \label{eq4.2}
\frac{3-5\beta}{\beta} (h|\nabla h|^2)^2 \phi
\le &~2 h^3 (h|\nabla h|^2) \left \langle \frac{(\beta-1) \nabla h}{\beta h} + \frac{\nabla \phi}{\phi}, \nabla \phi \right \rangle \nonumber \\
&~+h^3 (h|\nabla h|^2) (\partial_t - \Delta_f + 2\mathsf{k}) \phi + h^3 \mathsf{Z}_\Sigma \phi.  
\end{align}

We now proceed onto bounding from above each of the terms on the right-hand side of \eqref{eq4.2}. 
Again the argument proceeds by considering the two case $d(x_1, x_0,t_1) \le 1$ and $d(x_1, x_0,t_1) \ge 1$. 
However in view of certain similarities with the proof of Theorem \ref{thm1}, we shall remain brief, focusing 
on case two only and mainly on the differences. Towards this end the first two terms on the right-hand side 
of \eqref{eq4.2} are bounded directly in modulus by using the Cauchy-Schwarz followed by Young's inequality:  
\begin{align}
2 \frac{\beta-1}{\beta} h^2 (h|\nabla h|^2) \langle\nabla h, \nabla \phi\rangle 
&\le 2 \left| \frac{1-\beta}{\beta} h^2 (h|\nabla h|^2) \langle \nabla h, \nabla \phi \rangle \right| \nonumber \\
&\le \frac{2(1-\beta)}{\beta} \phi^{3/4} (h|\nabla h|^2)^{3/2} h^{3/2} \frac{|\nabla \phi|}{\phi^{3/4}} \nonumber \\
& \le \frac{1-\beta}{4\beta} (h|\nabla h|^2)^2 \phi + C(\beta) \frac{|\nabla \phi|^4}{\phi^3} h^6 \nonumber \\
&\le \frac{1-\beta}{4\beta} (h|\nabla h|^2)^2 \phi + \frac{C(\beta)}{R^4} h^6,   
\end{align}
where we have used $\sqrt W = \sqrt h |\nabla h|$ and in much the same way 
\begin{align}
2 \frac{|\nabla \phi|^2}{\phi}h^3 (h|\nabla h|^2) 
&\le 2 \phi^{1/2} (h|\nabla h|^2) \frac{|\nabla \phi|^2}{\phi^{3/2}} h^3 \nonumber \\
&\le \frac{1-\beta}{4\beta} (h|\nabla h|^2)^2 \phi + C(\beta) \frac{|\nabla \phi|^4}{\phi^3} h^6 \nonumber \\
&\le \frac{1-\beta}{4\beta} (h|\nabla h|^2)^2 \phi + \frac{C(\beta)}{R^4} h^6. 
\end{align}
For the $\mathsf{Z}_\Sigma$ term by noting 
$h^5 |\langle \nabla h, \Sigma_x \rangle| \le h^5 |\nabla h| |\Sigma_x| = h^{9/2} \sqrt h |\nabla h| |\Sigma_x|$ we have 
\begin{align}
h^3 \mathsf{Z}_\Sigma \phi &= \frac{2u \Sigma_u - (2-3\beta) \Sigma}{u} h^3 (h|\nabla h|^2) \phi 
+ \frac{2\beta \langle \nabla h, \Sigma_x \rangle}{u} h^5 \phi \nonumber \\
&~\le \left[ \frac{2u \Sigma_u - (2-3\beta) \Sigma}{u} \right]_+ h^3 (h|\nabla h|^2) \phi 
+ 2\beta h^{9/2} \sqrt{h|\nabla h|^2} \frac{|\Sigma_x|}{u} \phi \nonumber \\
&~\le \frac{1-\beta}{4 \beta} (h|\nabla h|^2)^2 \phi + C(\beta) [\mathsf{T}_\Sigma^2 + \mathsf{S}_\Sigma^{4/3}] h^6, 
\end{align}
where in the last line we have written $\mathsf{T}_\Sigma = \{ [2u \Sigma_u-(2-3\beta) \Sigma]/u\}_+$ 
and $\mathsf{S}_\Sigma = |\Sigma_x|/u$. For the terms involving $\Delta_f \phi$ and $\partial_t \phi$ 
we proceed similar to the proof of Theorem \ref{thm1}. Indeed recall from the Laplacian comparison 
theorem and the discussion on $(\partial_t-\Delta_f) \phi$ in Section \ref{sec2} that here we have 
the bounds  
\begin{align*}
- \Delta_f \phi &  \leq C \left( \frac{1}{R^2} + \frac{[\gamma_{\Delta_f}]_+}{R} + k_1 \right) \sqrt \phi \qquad 
\text{and } \qquad \partial_t \phi \le  C \left(\frac{1}{\tau} + k_2 \right) \sqrt \phi.
\end{align*}
Hence, by putting together all the three terms in this set in \eqref{eq4.2} we have, 
\begin{align} \label{Lap comp 2 sum}
h^4 |\nabla h|^2 \left( \partial_t - \Delta_f + 2\mathsf{k} \right) \phi 
& \le C h|\nabla h|^2 \sqrt \phi \left( \frac{1}{R^2} + \frac{[\gamma_{\Delta_f}]_+}{R} 
+ \frac{1}{\tau} + k_1 + k_2 + \mathsf{k} \right) h^3 \nonumber \\
& \le \frac{1-\beta}{4\beta} (h|\nabla h|^2)^2 \phi
+ C(\beta) \left( \frac{1}{R^4} + \frac{[\gamma_{\Delta_f}]_+^2}{R^2} + \frac{1}{\tau^2}+ k^2 \right) h^6.  
\end{align}

Having now estimated each of the individual terms on the right-hand side of \eqref{eq4.2} we proceed next by substituting these 
back into the inequality and finalising the estimate. To this end noting that the term $(3-5\beta) \phi (h|\nabla h|^2)^2/\beta$ on 
the left majorises, after summation, the term $(1-\beta) \phi (h|\nabla h|^2)^2/\beta$ on the right subject to $1-2\beta>0$, it 
follows from \eqref{eq4.2}-\eqref{Lap comp 2 sum}, upon rearranging terms and a basic set of calculations that  at the point 
$(x_1,t_1)$ we have   
\begin{align} \label{ing eq1}
(h|\nabla h|^2)^2 \phi \le C(\beta) &\left\{ \frac{1}{R^4} + \frac{[\gamma_{\Delta_f}]_+^2}{R^2} + \frac{1}{\tau^2} + k^2 
+ [\mathsf{T}_\Sigma^2 + \mathsf{S}_\Sigma^{4/3}] \right\} h^6.
\end{align}

By the maximality of $\phi h|\nabla h|^2$ at $(x_1, t_1)$ we have for any $x$ with $d(x, x_0,\tau) \le R/2$ the chain of inequalities 
\begin{equation}
(h^2|\nabla h|^4)(x,\tau) \le (\phi^2 h^2|\nabla h|^4)(x, \tau) \leq (\phi^2 h^2|\nabla h|^4)(x_1,t_1) \leq (\phi h^2|\nabla h|^4)(x_1,t_1),
\end{equation} 
[recall that $\phi(x, \tau) =1$ when $d(x,x_0; \tau) \leq R/2$]. Hence combining the latter with \eqref{ing eq1} and making note of the 
relations $h=u^\beta$ and $h |\nabla h|^2 = \beta^2 h^3 |\nabla u|^2/u^2$ this gives 
\begin{align*}
\frac{|\nabla u|^2}{u^{2-3\beta}} 
\le C(\beta) \Big( \sup_{Q_{R, T}} & u \Big)^{3\beta} \left\{ \frac{1}{R^2} 
+ \frac{[\gamma_{\Delta_f}]_+}{R} + \frac{1}{\tau} + k 
+ \sup_{Q_{R,T}} [\mathsf{T}_\Sigma + \mathsf{S}_\Sigma^{2/3}] \right\}. 
\end{align*}
The arbitrariness of $\tau \in (0, T]$ now gives the desired estimate for every $(x, t) \in Q_{R/2,T}$ with $t>0$. In particular setting 
$\beta =1/3$ and rearranging terms gives \eqref{eq1.12} when $q \equiv 0$ and $\mathsf{N} \equiv 0$. The passage to the general 
case with non-zero $q$ is analogous to Theorem \ref{thm1} upon noting 
$\mathsf{T}_{\Sigma+qu} \le \mathsf{T}_{\Sigma} + \mathsf{T}_{qu} = \mathsf{T}_{\Sigma} + q_+$ and 
$\mathsf{S}_{\Sigma+qu} \le \mathsf{S}_{\Sigma} + \mathsf{S}_{qu} = \mathsf{S}_{\Sigma} + |\nabla q|$. 
The rest of the proof is similar and thus abbreviated. \qed

\section{Proof of the parabolic Li-Yau type estimate in Theorem \ref{thm28}}  \label{sec6}

This section is devoted to the proof of the nonlinear version of the Li-Yau estimate (also known as the differential Harnack 
inequality) in Theorem \ref{thm28}. As the proof is quite involved and requires several intermediate steps, for the sake of 
clarity and convenience, we break this into three subsections, focusing first on deriving and establishing some of the 
necessary tools and identities and then finalising the proof in the last subsection.

\subsection{Some basic identities on evolutionary metrics and potentials} The results and identities proved 
here will be repeatedly used throughout the section. Our first task is to obtain a relationship between 
$\partial_t \Delta_f h$ and $\Delta_f \partial_t h$ for a smooth function $h=h(x,t)$ 
given that both the metric $g$ and the  potential $f$ are time dependent. The following lemma gives the 
required ingredients. For convenience in writing we hereafter set 
\begin {align} \label{2.4}
\frac{\partial g}{\partial t} (x,t)= 2v (x,t). 
\end{align}

\begin{lemma} \label{geometric-evolution-lemma-one}
With the notation introduced above for every smooth function $h=h(x,t)$ we have the relations 
\begin{equation} \label{norm-grad-evolve-equation}
\partial_t |\nabla h|^2 = -2v(\nabla h, \nabla h) + 2 \langle \nabla h, \nabla \partial_t h \rangle, 
\end{equation}
\begin{equation} \label{Lap-evolve-equation}
\partial_t \Delta h 
= \Delta \partial_t h - 2 \langle v, \nabla^2 h \rangle 
- \langle 2{\rm div}\,v - \nabla ({\rm Tr}_g\,v), \nabla h \rangle, 
\end{equation}
where $({\rm div}\,v)_k=g^{ij} \nabla_i v_{jk}$ and ${\rm Tr}_g\,v = g^{ij} v_{ij}$. 
\end{lemma}

\begin{proof}
The identity in \eqref{norm-grad-evolve-equation} follows by first writing $|\nabla h|^2=g^{ij} \nabla_i h \nabla_j h$ and then taking $\partial_t$ 
making note of $\partial_t g^{ij} = - 2 g^{ik} g^{j \ell} v_{k \ell} = -2v^{ij}$. For the identity in \eqref{Lap-evolve-equation} we first recall the relation  
\begin{align} \label{dt-Ch-symb}
\frac{\partial}{\partial t} \Gamma^k_{ij} = g^{k \ell} (\nabla_i v_{j \ell} + \nabla_j v_{i \ell} - \nabla_\ell v_{ij}).
\end{align}
Then direct differentiation and making use of \eqref{dt-Ch-symb} leads to  
\begin{align}
\partial_t \Delta h &= \partial_t [g^{ij} (\nabla_i \nabla_j - \Gamma^k_{ij} \nabla_k )h] \nonumber \\
&= (\partial_t g^{ij}) (\nabla_i \nabla_j - \Gamma^k_{ij} \nabla_k) h 
+  g^{ij} (\nabla_i \nabla_j - \Gamma^k_{ij} \nabla_k ) \partial_t h - g^{ij} (\partial_t \Gamma^k_{ij}) \nabla_k h \nonumber \\
&= -2v^{ij} (\nabla_i \nabla_j - \Gamma^k_{ij} \nabla_k) h 
+  g^{ij} (\nabla_i \nabla_j - \Gamma^k_{ij} \nabla_k ) \partial_t h - g^{ij} (\partial_t \Gamma^k_{ij}) \nabla_k h \nonumber \\
&= -2 \langle v, \nabla^2 h \rangle + \Delta \partial_t h - g^{ij} g^{k \ell} [\nabla_i v_{j \ell} + \nabla_j v_{i \ell} - \nabla_\ell v_{ij}] \nabla_k h \nonumber \\
&= \Delta \partial_t h - 2 \langle v, \nabla^2 h \rangle - g^{k \ell} [2 g^{ij} \nabla_i v_{j \ell} - \nabla_\ell ({\rm Tr}_g\,v)] \nabla_k h \nonumber \\
&= \Delta \partial_t h - 2 \langle v, \nabla^2 h \rangle - \langle 2 {\rm div}\,v - \nabla ({\rm Tr}_g\,v), \nabla h \rangle
\end{align}
which is the desired identity. 
\end{proof}

\begin{lemma} \label{geometric-evolution-lemma-two}
Subject to the notation in $\eqref{2.4}$, for every pair of smooth functions $f=f(x,t)$ and $h=h(x,t)$ we have 
\begin{equation} \label{inner-product-grad-evolve-equation}
\partial_t \langle \nabla f, \nabla h \rangle 
= -2v(\nabla f, \nabla h) + \langle \nabla \partial_t f, \nabla h \rangle + \langle \nabla f, \nabla \partial_t h \rangle, 
\end{equation}
\begin {align} \label{Lap-f-evolve-equation}
\partial_t( \Delta_f h ) 
=&~\Delta_f (\partial_t h)  -2 \langle v , \nabla ^2 h \rangle - \langle 2 {\rm div} v - \nabla ({\rm Tr}_g v),\nabla h \rangle \nonumber \\
&- \langle \nabla \partial_t f , \nabla h \rangle + 2v (\nabla f , \nabla h). 
\end{align}
\end{lemma}

\begin{proof}
The first identity follows similar to \eqref{norm-grad-evolve-equation} in Lemma \ref{geometric-evolution-lemma-one}. 
For the second identity we proceed by directly calculating $\partial_t \Delta_f h$ whilst making note of 
\eqref{Lap-evolve-equation} and \eqref{inner-product-grad-evolve-equation}
\begin {align} 
\partial_t( \Delta_f h) =&~\partial_t(\Delta h - \langle \nabla f , \nabla h \rangle) \nonumber \\
=&~\Delta (\partial_t h)  -2 \langle v , \nabla ^2 h \rangle - \langle 2 {\rm div} v - \nabla ({\rm Tr}_g v),\nabla h \rangle \nonumber \\
&- \langle \nabla \partial_t f , \nabla h \rangle - \langle \nabla f, \nabla \partial_t h \rangle + 2v \langle \nabla f , \nabla h \rangle \nonumber \\
=&~\Delta_f (\partial_t h)  -2 \langle v , \nabla ^2 h \rangle - \langle 2 {\rm div} v - \nabla ({\rm Tr}_g v),\nabla h \rangle \nonumber \\
&- \langle \nabla \partial_t f , \nabla h \rangle + 2v (\nabla f , \nabla h). 
\end{align}
The proof is thus complete.
\end{proof}

\begin{lemma} \label{LiYau-basiclemma}
Suppose $a,b,z \in \mathbb{R}$, $c, y>0$ and $\lambda>1$ are arbitrary constants such that $y-\lambda z >0$. Then for any 
$\varepsilon \in (0,1)$ we have
\begin{align} \label{eq6.9} 
(y-z)^2& - a \sqrt y (y-\lambda z) - b y - c \sqrt y \nonumber\\
&\ge (y-\lambda z)^2/\lambda^2-a^2\lambda^2 (y-\lambda z)/[8 (\lambda-1)] \nonumber \\
& - (\lambda^2 b^2)/[4(1-\varepsilon)(\lambda-1)^2]
-(3/4) c^{4/3} [\lambda^2/(4 \varepsilon (\lambda -1)^2)]^{1/3}.  
\end{align}
\end{lemma}

\begin{proof}
Starting from the expression on the left-hand side in \eqref{eq6.9} we can write for any $\delta, \varepsilon$ by basic considerations
\begin {align} \label{eq6.10}
(y-z)^2 &- a \sqrt y (y-\lambda z) - b y - c \sqrt y\nonumber\\
=&~(1-\varepsilon-\delta)y^2-(2-\varepsilon \lambda)yz+z^2+(\varepsilon y - a \sqrt y)(y-\lambda z) 
+ \delta y^2 - by - c \sqrt y \nonumber\\
=&~(1/\lambda - \varepsilon/2)(y-\lambda z)^2+(1-\varepsilon-\delta-1/\lambda+\varepsilon/2)y^2
+(1- \lambda + \varepsilon \lambda^2/2)z^2\nonumber\\
&+ (\varepsilon y - a \sqrt y)(y-\lambda z) + \delta y^2 - by - c \sqrt y. 
\end{align}   
In particular setting $\delta = (1/\lambda-1)^2$ and $\varepsilon = 2-2/\lambda-2(1/\lambda-1)^2 = 2(\lambda-1)/\lambda^2$ 
gives $1-\varepsilon-\delta-1/\lambda+\varepsilon/2=0$ and $1-\lambda + \varepsilon\lambda^2/2=0$ and so by making note 
of the inequality $\varepsilon y - a \sqrt y \ge - a^2/(4 \varepsilon)$ with $\varepsilon=2(\lambda-1)/\lambda^2>0$ we can 
deduce from 
\eqref{eq6.10} that
\begin {align}\label{8.36}   
(y-z)^2& - a \sqrt y (y-\lambda z) - b y - c \sqrt y \\
&\ge (y-\lambda z)^2/\lambda^2-a^2\lambda^2 (y-\lambda z)/[8 (\lambda-1)] 
+ (\lambda -1)^2 y^2/\lambda^2 - by- c \sqrt y. \nonumber 
\end{align} 
Next, considering the last three terms in the above inequality only, we can write, for any $\varepsilon \in (0,1)$,
\begin {align} 
(\lambda -1)^2 &y^2/\lambda^2 - by - c \sqrt y \nonumber\\
&=(\lambda -1)^2 y^2/\lambda^2 - (1-\varepsilon)(\lambda-1)^2 y^2/\lambda^2 
+ (1-\varepsilon)(\lambda-1)^2 y^2/\lambda^2- by - c \sqrt y \nonumber \\
&\ge (\lambda -1)^2y^2/\lambda^2 - (1-\varepsilon)(\lambda-1)^2 y^2/\lambda^2 - 
(\lambda^2 b^2)/[4(1-\varepsilon)(\lambda-1)^2] - c \sqrt y\nonumber\\
&\ge \varepsilon (\lambda -1)^2y^2/\lambda^2 - 
(\lambda^2 b^2)/[4(1-\varepsilon)(\lambda-1)^2] - c \sqrt y\nonumber\\
&\ge -(3/4) c^{4/3} [\lambda^2/(4 \varepsilon (\lambda -1)^2)]^{1/3} 
- (\lambda^2 b^2)/[4(1-\varepsilon)(\lambda-1)^2] 
\end{align} 
where above we have made use of $(1-\varepsilon)(\lambda-1)^2 y^2/\lambda^2 - by \ge -(\lambda^2 b^2)/[4(1-\varepsilon)(\lambda-1)^2]$ 
and $\varepsilon (\lambda -1)^2y^2/\lambda^2 - c \sqrt y \ge  -(3/4) c^{4/3} [\lambda^2/(4 \varepsilon (\lambda -1)^2)]^{1/3}$  
to deduce the first and last inequalities respectively. Substituting back in \eqref{8.36} gives the desired inequality.  
\end{proof}

\subsection{Evolution of a Harnack quantity}

In this subsection we introduce a Harnack quantity built out of the solution $u$ and establish a parabolic inequality 
by considering its evolution under the weighted heat operator.

\begin{lemma} \label{Lemma 2.1}
Let $u$ be a positive solution to the equation $(\partial_t - \Delta_f) u = \Sigma (t,x,u)$ and let $G=G(x,t)$ be defined by 
\begin {align} \label{2.3}
G= t [|\nabla h|^2 - \lambda \partial_t h + \lambda e^{-h} \Sigma(t,x,e^h)], \qquad t \ge 0, 
\end{align}  
where $h =\log u$ and $\lambda >1$ is a fixed constant. Suppose that the metric-potential pair $(g, f)$ is time dependent 
and of class $\mathscr{C}^2$ and that we have $\eqref{2.4}$. Then $G$ satisfies the evolution equation 
\begin{align} \label{9.4a}
\square G = (\partial_t - \Delta_f) G=& ~2 \langle \nabla h , \nabla G \rangle 
- 2 t|\nabla ^2 h |^2 - 2t {\mathscr Ric}_f^m (\nabla h, \nabla h)  \nonumber\\ 
& + 2t (\lambda -1) v (\nabla h, \nabla h) 
- 2t \langle \nabla f , \nabla h \rangle^2/(m-n) + G/t \nonumber \\
& + 2 \lambda t [\langle v , \nabla ^2 h \rangle 
+ \langle {\rm div} v - (1/2) \nabla ({\rm Tr}_g v),\nabla h \rangle] \\
& + \lambda t[\langle \nabla \partial_t f , \nabla h \rangle- 2v (\nabla f , \nabla h) ] \nonumber \\
& -2t (\lambda -1) \langle \nabla h, \nabla [e^{-h}\Sigma(t,x,e^h)] \rangle 
- \lambda t \Delta_f [e^{-h} \Sigma (t,x,e^h)]. \nonumber 
\end{align}
\end{lemma}

\begin{proof}
Referring to the equation for $u$ an easy calculation shows that $h$ in turn satisfies the equation 
\begin {align} \label{2.2}
\square h = (\partial_t - \Delta_f) h = |\nabla h|^2 + e^{-h} \Sigma(t,x,e^h).
\end{align} 
Moreover, using \eqref{2.3} and \eqref{2.2} it is easily seen that the following relation emerges between $G$, $|\nabla h|^2$ and $\Delta_f h$: 
\begin {align} \label{9.6a}
\Delta_f h &= - \left[ |\nabla h|^2/\lambda - \partial_t h + e^{-h} \Sigma(t,x,e^h) \right] 
- (\lambda-1)|\nabla h|^2/\lambda \nonumber \\
&= - G/(\lambda t) - (\lambda-1)|\nabla h|^2/\lambda, \qquad t>0. 
\end{align}

Now having these identities and relations in place we next proceed onto applying the weighted heat operator $\partial_t - \Delta_f$ 
to the Harnack quantity $G$ given by \eqref{2.3}. Towards this end we first note that     
\begin {align} \label{Witten-Lap-F-DHI-equation}
\Delta_f G = t (\Delta_f |\nabla h|^2 - \lambda \Delta_f (\partial_t h) + \lambda \Delta_f [e^{-h}\Sigma(t,x,e^h)]).
\end{align}

As for the first term on the right by recalling the weighted Bocnher-Weitzenb\"ock formula as applied to $h$ we have 
\begin {align} 
\Delta_f |\nabla h|^2/2 = |\nabla^2 h|^2 + \langle \nabla h , \nabla \Delta_f h \rangle 
+ {\mathscr Ric}_f^m (\nabla h, \nabla h) + \langle \nabla f , \nabla h \rangle^2/(m-n),
\end{align}
and so upon substituting back in \eqref{Witten-Lap-F-DHI-equation} and making use of \eqref{Lap-f-evolve-equation} this gives
\begin {align} \label{2.11}
\Delta_f G =&~t \left[ 2 |\nabla ^2 h |^2 +2 \langle \nabla h, \nabla \Delta_f h \rangle 
+2 {\mathscr Ric}_f^m (\nabla h, \nabla h)+2 \frac{\langle \nabla f , \nabla h \rangle^2}{m-n} \right]\nonumber\\
& -\lambda t \partial_t (\Delta_f h) -2\lambda t \left[ \langle v , \nabla ^2 h \rangle + \langle {\rm div} v 
- \frac{1}{2}\nabla ({\rm Tr}_g v),\nabla h \rangle \right] \nonumber \\
& -\lambda t [\langle \nabla \partial_t f , \nabla h \rangle - 2v ( \nabla f , \nabla h )] 
+ \lambda t \Delta_f [e^{-h} \Sigma (t,x,e^h)].
\end{align}

Now referring to the sum on the right the contributions of the second and fifth terms, modulo a factor $t$ and upon 
using \eqref{2.2} and \eqref{9.6a} can be simplified and re-written as, 
\begin {align} 
2 \langle \nabla h, \nabla \Delta_f h \rangle -& \lambda \partial_t(\Delta_f h)  \nonumber\\
= &~2 \langle \nabla h, \nabla \Delta_f h \rangle 
- 2 (\lambda -1) v (\nabla h , \nabla h) \nonumber\\
&- [-(t\partial_t G - G)/t^2 
- 2 (\lambda -1)\langle \nabla h , \nabla (\partial_t h)\rangle] \nonumber \\
= &~2 (\lambda -1) \langle \nabla h , \nabla [\Delta_f h +|\nabla h|^2 
+ e^{-h}\Sigma (t,x,e^h)]\rangle \nonumber\\
&+2 \langle \nabla h, \nabla \Delta_f h \rangle + (t \partial_t G - G)/t^2 - 2 (\lambda -1) v (\nabla h , \nabla h) \nonumber \\
= &~2 \lambda \langle \nabla h , \nabla [ -G/(\lambda t) - (\lambda -1)|\nabla h|^2/\lambda] \rangle \nonumber \\
&+ 2(\lambda -1) \langle \nabla h, \nabla |\nabla h|^2 \rangle 
+ (t \partial_t G - G)/t^2 \nonumber \\
&+2 (\lambda -1) \langle \nabla h , \nabla [e^{-h}\Sigma (t,x,e^h)] \rangle
- 2 (\lambda -1) v (\nabla h , \nabla h)\nonumber\\
= &~2 (\lambda -1) [\langle \nabla h , \nabla [e^{-h}\Sigma (t,x,e^h)] \rangle 
- v (\nabla h , \nabla h)] \nonumber \\
&+ (t\partial_t G - G)/t^2 - 2 \langle \nabla h , \nabla G \rangle/t. 
\end{align}
Therefore substituting this back into \eqref {2.11} gives 
\begin {align}\label{1.13} 
\Delta_f G = &~\partial_t G + 2 t|\nabla ^2 h |^2 
+2t (\lambda -1)[ \langle \nabla h, \nabla[e^{-h}\Sigma(t,x,e^h)] \rangle - v (\nabla h, \nabla h)] \nonumber\\ 
& - G/t -2 \langle \nabla h , \nabla G \rangle 
+ 2t {\mathscr Ric}_f^m (\nabla h, \nabla h) +2t \langle \nabla f , \nabla h \rangle^2/(m-n)  \nonumber\\
&- 2\lambda t \langle v , \nabla ^2 h \rangle 
- \lambda t \langle 2 {\rm div} v - \nabla ({\rm Tr}_g v),\nabla h \rangle 
- \lambda t \langle \nabla \partial_t f , \nabla h \rangle\nonumber\\
&+2 \lambda tv (\nabla f , \nabla h) + \lambda t \Delta_f[ e^{-h} \Sigma (t,x,e^h)]
\end{align}
which upon a rearrangement of term leads to the desired conclusion. 
\end{proof}

\begin{lemma} \label{estimate-second-type-lemma}
Let $u$ be a positive solution to $(\partial_t - \Delta_f) u = \Sigma (t,x,u)$ and let $G=G(x,t)$ be as in \eqref{2.3}. 
Suppose that the metric-potential pair $(g, f)$ is time dependent and of class $\mathscr{C}^2$. Assume additionally 
that ${\mathscr Ric}_f^m (g) \ge -(m-1) k_1 g$ and we have $[$cf. \eqref{2.4}$]$
\begin{equation} \label{eq-v-6.22}
- \underline k_2 g \le v \le \overline k_2 g, \qquad |\nabla v| \le k_3
\end{equation} 
for suitable $k_1, \underline k_2, \overline k_2$ and $k_3 \ge 0$. Then
\begin {align}\label{1.20}
\square G = (\partial_t - \Delta_f) G \le &-t (\Delta_f h)^2/m + G/t + 2 \langle \nabla h, \nabla G \rangle \nonumber \\
& + 2t [(m-1)k_1 +(\lambda -1) \overline k_2] |\nabla h|^2\nonumber\\
&+ \lambda^2 n t(\underline k_2 + \overline k_2)^2 + 3 \lambda t\sqrt n k_3 |\nabla h| \nonumber \\
& + \lambda t \langle \nabla \partial_t f , \nabla h \rangle - 2\lambda t  v ( \nabla f , \nabla h ) \nonumber\\
&- 2 (\lambda -1) t \langle \nabla h,\nabla [e^{-h} \Sigma(t,x,e^h)] \rangle \nonumber \\
& - \lambda t\Delta_f [e^{-h} \Sigma(t,x,e^h)].
\end{align}
\end{lemma}

\begin{proof}
Firstly from \eqref{eq-v-6.22} we deduce that $|v|^2  \le (\underline k_2 + \overline k_2)^2 |g|^2 = n(\underline k_2 + \overline k_2)^2$ 
and so
\begin {align} 
|\lambda \langle v , \nabla^2 h \rangle| \leq \frac{1}{2} |\nabla ^2 h|^2 + \frac{1}{2}\lambda ^2 |v|^2 \leq \frac{1}{2} |\nabla ^2 h|^2 
+\frac{1}{2} \lambda ^2 n(\underline k_2 + \overline k_2)^2.
\end{align} 
Next in view of $|g^{ij} (2 \nabla_i v_{j \ell} - \nabla_\ell v_{ij})| 
\le 3|g| |\nabla v|$, $|g|=\sqrt n$ and $|\nabla v| \le k_3$ we can write 
\begin {align*} 
\left|{\rm div} v - \frac{1}{2} \nabla ({\rm Tr}_g v)\right| &= \left| g^{ij} \nabla_i v_{j \ell} - \frac{1}{2} g^{ij} \nabla_\ell v_{ij} \right| 
= \frac{1}{2} \left| g^{ij} (2 \nabla_i v_{j \ell} - \nabla_\ell v_{ij}) \right| \le \frac{3}{2} \sqrt n k_3. 
\end{align*} 
The conclusion now follows at once by recalling \eqref{9.4a} in Lemma \ref{Lemma 2.1}, making note of  
\begin{align}
|\nabla^2 h|^2+\frac{\langle \nabla f , \nabla h \rangle^2}{m-n} 
\ge \frac{(\Delta h)^2}{n}+\frac{\langle \nabla f , \nabla h \rangle^2}{m-n} 
\ge\frac{(\Delta_f h)^2}{m},  
\end{align} 
and the Bakry-\'Emery curvature lower bound ${\mathscr Ric}_f^m(g) \ge -(m-1) k_1 g$ in the lemma. 
 \end{proof}

\subsection{Proof of the local estimate in Theorem \ref{thm28}}

Having all the ingredients and necessary tools at our disposal we now come to the proof of the main estimate. 
The proof of the theorem is broken into two parts, first, the case $q \equiv 0$, and then the general case that 
follows from the latter by incorporating the $qu$ term in the nonlinearity $\Sigma$ and then re-calculating the 
resulting bounds accordingly. The idea for the first part is 
to combine the inequality established in Lemma \ref{estimate-second-type-lemma} along with a 
localisation argument by utilising a suitable cut-off function. The estimates on the cut-off function in turn makes 
use of the generalised Laplacian comparison theorem and the Bakry-Eme\'ry generalised Ricci curvature lower 
bound as will be described in detail in the course of the proof. 
To this end we first set $q \equiv 0$. We pick a reference point $x_0 \in M$ and fix $R, T>0$ and $0<\tau \le T$. 
As before we denote by $\varrho(x,t)=d(x,x_0,t)$ the geodesic radial variable at time $t$ in reference to 
$x_0$. For the sake of localisation we consider first a function $\bar{\psi} = \bar{\psi}(s)$ 
on the half-line $s \ge 0$ ({\it see} Lemma \ref{psi lemma} below) and then for $x \in M$ and $t \ge 0$ set 
\begin{align}\label{9.31}
\psi(x,t)=\bar{\psi} \left( \frac{\varrho (x, t)}{R} \right). 
\end{align}
The existence of $\bar{\psi}=\bar \psi(s)$ as 
used in \eqref{9.31} and its properties is granted by the following straightforward and standard statement.

\begin{lemma} \label{psi lemma} There exists a function $\bar{\psi}:[0,\infty) \to \mathbb{R}$ 
verifying the following properties:  
\begin{enumerate}[label=$(\roman*)$]
\item $\bar\psi$ is of class $\mathscr{C}^2 [0, \infty)$, 
\item $0 \le \bar\psi(s) \le 1$ for $0 \le s < \infty$ with $\bar\psi \equiv 1$ for $s \le 1$ and $\bar\psi \equiv 0$ for $s \ge 2$, 
\item $\bar\psi$ is non-increasing $($specifically, $\bar\psi' \le 0$$)$ and additionally, for suitable constants $c_1, c_2>0$, 
satisfies the bounds 
\begin{align} \label{9.30}
- c_1 \le \frac{\bar{\psi}^{'}}{\sqrt{ \bar{\psi}}} \le 0, \qquad and \qquad  \bar{\psi}^{''} \ge -c_2, 
\end{align}
on the half-line $[0, \infty)$.
\end{enumerate}
\end{lemma}

It is evident from $(ii)$ that $\psi \equiv 1$ for when $0 \le \varrho(x, t) \le R$ and $\psi \equiv 0$ for 
when $\varrho(x,t) \ge 2R$. 
Let us now consider the spatially localised function $\psi G$ where $G$ is as in \eqref{2.3}. 
We denote by $(x_1,t_1)$ the point where this function attains its maximum over the compact 
set $\{d(x,x_0,t) \le 2R, 0 \le t \le \tau\}$. We will also assume that $[\psi G] (x_1, t_1)>0$ as otherwise 
the desired estimate is trivially true as a result of $G \le 0$. It thus follows that $t_1>0$ and 
$d(x_1, x_0, t_1) < 2R$ and so at the maximum point $(x_1, t_1)$ we have the relations 
\begin{align} \label{9.32}
\nabla (\psi G)=0, \qquad \partial_t (\psi G) \ge 0, \qquad \Delta(\psi G) \le 0, \qquad \Delta_f (\psi G) \le 0.
\end{align}
Starting with the basic identity 
$\Delta_f (\psi G) = G \Delta_f \psi + 2 \langle \nabla \psi, \nabla G \rangle + \psi \Delta_f G$ 
and making note of the relations \eqref{9.32} at the maximum point $(x_1, t_1)$ we can write 
\begin{align} \label{9.34}
0 &\ge \Delta_f (\psi G) 
= G \Delta_f \psi + (2/\psi) \langle \nabla \psi , \nabla (\psi G) \rangle - 2(|\nabla \psi |^2/\psi) G + \psi \Delta_f G \nonumber \\ 
& \ge G \Delta_f \psi -2(|\nabla \psi |^2/\psi) G + \psi \Delta_f G.
\end{align}
Now from \eqref{9.31} we deduce $\nabla \psi = (\bar\psi'/R) \nabla \varrho$ and 
$\Delta \psi = \bar\psi'' |\nabla \varrho|^2/R^2+\bar\psi' \Delta \varrho/R$ and so 
\begin{equation} \label{9.35}
\Delta_f \psi = \Delta \psi - \langle \nabla f, \nabla \psi \rangle = (\bar\psi''/R^2) |\nabla \varrho|^2 + (\bar\psi'/R) \Delta_f \varrho.
\end{equation} 
Since ${\mathscr Ric}_f^m \geq -(m-1) k_1 g$ we have $\Delta_f \varrho \le (m-1) \sqrt {k_1}\coth (\sqrt{k_1} \varrho)$ 
({\it cf.} \cite{[WeW09]}) and so from \eqref{9.35} we have:
\begin{align}
\Delta_f \psi \ge (\bar\psi''/R^2) + (m-1)(\bar\psi'/R) \sqrt{k_1} \coth (\sqrt{k_1} \varrho).
\end{align} 
Next $\coth (\sqrt{k_1} \varrho) \le \coth (\sqrt{k_1} R)$ and $\sqrt{k_1} \coth (\sqrt{k_1} R) \le (1 + \sqrt{k_1} R)/R$ 
for $R \le \varrho \le 2R$ and therefore 
$(m-1) \bar\psi' \sqrt{k_1} \coth (\sqrt{k_1}\varrho)\ge (m-1) [1+\sqrt{k_1}R] (\bar\psi'/R)$. Hence
\begin{align} \label{9.38}
\Delta_f \psi &\ge \frac {1}{R^2}\bar{\psi}^{''}+\frac{(m-1)}{R} \left( \frac{1}{R}+\sqrt{k_1} \right) \bar{\psi}^{'} \nonumber \\
&\ge -\frac{1}{R^2}[c_2 +(m-1) c_1(1+R \sqrt{k_1})].
\end{align}
Thus returning to \eqref{9.34}, invoking \eqref{1.20} and making note of \eqref{9.38}, we obtain, at the maximum point $(x_1, t_1)$, the inequality 
\begin{align}\label {2.25}
0 \ge&~\Delta_f (\psi G) = G \Delta_f \psi -2 (|\nabla \psi |^2/\psi) G + \psi \Delta_f G \nonumber \\
\ge&-(1/R^2) [c_2 +(m-1) c_1(1+R \sqrt{k_1})] G -2(|\nabla \psi |^2/\psi) G + \psi \partial_t G \nonumber\\
&+ \psi \bigg[ (t_1/m) (\Delta_f h)^2 - (G/t_1) -2 \langle \nabla h, \nabla G \rangle -2 [(m-1)k_1 
+(\lambda -1) \overline k_2] t_1 |\nabla h|^2\nonumber\\
&-\lambda^2 n t_1(\underline k_2 + \overline k_2)^2 -3 \sqrt n k_3 \lambda t_1 |\nabla h| -\lambda t_1 \langle \nabla \partial_t f , \nabla h \rangle 
+ 2\lambda t_1  v ( \nabla f , \nabla h)\nonumber\\
&+ 2 (\lambda -1) t_1 \langle \nabla h,\nabla (e^{-h} \Sigma)\rangle +\lambda t_1 \Delta_f (e^{-h} \Sigma) \bigg],
\end{align}
where for the sake of convenience we have abbreviated the arguments of $\Sigma=\Sigma(t,x,e^h)$. We now proceed by bounding the 
individual terms in the last inequality. Here, starting from the last term on the first line, upon recalling \eqref{9.31}, we have, 
\begin{align}
\psi \partial_t G =  \partial_t(\psi G) - G \partial_t \psi = \partial_t(\psi G) - G \bar\psi' \left( \frac{\varrho}{R} \right) \frac{\partial_t \varrho}{R},     
\end{align}
and so recalling that at the maximum point $(x_1, t_1)$ we have $\partial_t (\psi G) \ge 0$, by restricting to this point the latter results in, 
\begin{align} \label{9.41}
\psi \partial_t G &\ge - G \partial_t \psi = - G \bar{\psi}^{'} \left( \frac{\varrho}{R} \right) \frac{\partial_t \varrho}{R} 
\ge G \bar{\psi}^{'} \left( \frac{\varrho}{R} \right) \frac{\underline k_2 \varrho}{R} \nonumber \\
& \ge -c_1 \underline k_2 \sqrt{\bar\psi \left( \frac{\varrho}{R} \right)} \frac{\varrho}{R} G \ge - c_1 \underline k_2 G. 
\end{align}

Here we point out that we have used $(iii)$ in the set of assumptions on $\bar\psi$, specifically, the first inequality in \eqref{9.30}. 
Note also that in deducing the second inequality on the first line we have made use of the relation 
\begin{align}
\frac{\partial}{\partial t} \varrho(x,t) &= \frac{\partial}{\partial t} \int_0^1 |\gamma'|_{g_t} \, ds \nonumber \\
&= \frac{\partial}{\partial t} \int_0^1 \sqrt{g_t(\gamma', \gamma')} \, ds 
= \frac{1}{2} \int_0^1 \frac{[\partial_t g_t](\gamma', \gamma')}{\sqrt{g_t(\gamma', \gamma')}} \, ds \nonumber \\
&= \int_0^1 \frac{v(\gamma', \gamma')}{|\gamma'|_{g_t}} \, ds \ge - \underline k_2 \int_0^1 |\gamma'|_{g_t}\, ds 
= - \underline k_2 \varrho(x,t),
\end{align}
where $\gamma=\gamma (s)$ with $0 \le s \le 1$ is a geodesic curve with respect to $g_t$ at fixed $t$ joining 
$\bar x$ to $x_1$, that is, $\gamma (0)= \bar x$ and $\gamma (1)= x_1$. Lastly we have used 
$v \ge - \underline k_2 g$ to obtain the final inequality.

Referring again to the last inequality in \eqref{2.25} and bounding the individual terms it is evident that 
$|\nabla \psi |^2/\psi = (\bar\psi'^2/\bar\psi) (|\nabla \varrho|^2/R^2) \le  c_1^2/R^2$,  
where we have again made use of $(iii)$ in the set of assumptions on $\bar\psi$. As a result using the above in \eqref{2.25} 
and rearranging terms we have 
\begin{align} 
0 \ge & -\frac{1}{R^2}[c_2 +(m-1) c_1(1+R \sqrt{k_1})] G -2 \frac{c_1^2}{R^2} G - c_1 \underline k_2 G  \nonumber\\
& + \frac{t_1}{m} \psi (\Delta_f h)^2 - \psi \frac{G}{t_1} - 2\psi \langle \nabla h, \nabla G \rangle 
- 2t_1 [(m-1)k_1 +(\lambda -1) \overline k_2] \psi |\nabla h|^2 \nonumber\\
& - \psi t_1 \bigg[\lambda^2 n (\underline k_2 + \overline k_2)^2 + 3 \sqrt n k_3 \lambda |\nabla h|\bigg] 
+ \lambda \psi t_1 \bigg[2 v ( \nabla f , \nabla h) - \langle \nabla \partial_t f , \nabla h \rangle \bigg] \nonumber\\
& + \psi t_1 \bigg[ 2 (\lambda -1) \langle \nabla h,\nabla (e^{-h} \Sigma) \rangle +\lambda \Delta_f (e^{-h} \Sigma) \bigg].
\end{align}
Using $\psi \langle \nabla h, \nabla G \rangle = -G \langle \nabla h , \nabla \psi \rangle 
\le G |\nabla h| |\nabla \psi| \le c_1 (\sqrt \psi/R) G |\nabla h|$ as $\nabla (\psi G) =0$ at $(x_1, t_1)$
and using $v (\nabla f , \nabla h) \ge - \underline k_2 |\nabla f| |\nabla h|$ 
and $\langle \nabla \partial_t f , \nabla h \rangle \le |\nabla \partial_t f| |\nabla h|$ 
together with $3 k_3 \sqrt n \lambda |\nabla h| \le 2 n k_3 \lambda ^2 + 2 k_3 |\nabla h|^2$ 
it then follows that  
\begin{align} \label{9.45}
0 \ge &- [c_2 +(m-1) c_1(1+R \sqrt{k_1})+2c_1^2] G/R^2 \nonumber \\
&- c_1 \underline k_2 G - 2c_1 (\sqrt \psi/R) G |\nabla h| + t_1 (\psi/m) (\Delta_f h)^2 \nonumber \\
&- \psi G/t_1 - 2t_1 [(m-1)k_1 +(\lambda -1) \overline k_2] \psi |\nabla h|^2 \nonumber\\
&- t_1 \psi [\lambda^2 n (\underline k_2 + \overline k_2)^2 + 2 n k_3 \lambda ^2 + 2 k_3 |\nabla h|^2] \nonumber \\
&- \lambda t_1 \psi [2 \underline k_2 |\nabla f| |\nabla h| + |\nabla \partial_t f| |\nabla h|] \nonumber\\
&+ t_1 \psi [2 (\lambda -1) \langle \nabla h,\nabla (e^{-h} \Sigma) \rangle +\lambda \Delta_f (e^{-h} \Sigma)].
\end{align}

Next multiplying \eqref{9.45} through by $t_1 \psi$, making note of \eqref{2.2} and 
rearranging terms gives 
\begin{align}\label{1.37}
0 \ge&-t_1 \psi G ([c_2 +(m-1) c_1(1+R \sqrt{k_1})] + 2c_1^2)/R^2 - \psi^2 G \nonumber\\
&- c_1 \underline k_2 t_1 \psi G + t_1^2 (\psi ^2/m) \left[|\nabla h|^2 + e^{-h} \Sigma 
- \partial_t h\right]^2 \nonumber \\ 
&- 2c_1 t_1 \psi (\sqrt \psi/R) G|\nabla h| 
-2 t_1^2 \left[(m-1)k_1 +(\lambda -1) \overline k_2 + k_3 \right] \psi ^2 |\nabla h|^2 \nonumber\\
&- n \lambda^2 t_1^2 \psi ^2[ (\underline k_2 + \overline k_2)^2+ 2 k_3] 
-\lambda t_1^2 \psi^2 |\nabla \partial_t f| |\nabla h| 
- 2 \lambda t_1^2 \psi ^2 \underline k_2 |\nabla f| |\nabla h| \nonumber \\
&+ t_1^2 \psi ^2 [2(\lambda -1)\langle \nabla h,\nabla (e^{-h} \Sigma) \rangle 
+ \lambda \Delta_f (e^{-h} \Sigma)].
\end{align}

We now go through some calculations relating to the nonlinearity $\Sigma=\Sigma(t,x,u)$ with $u=e^h$ 
and $h=h(x,t)$, abbreviating the arguments $(t, x, u)$ for convenience. Firstly, it is seen, by calculating 
in local coordinates or directly, that 
\begin {align} \label{eq6.44}
\nabla \Sigma = \Sigma_x + e^h \Sigma_u \nabla h, \qquad \Sigma_x=(\Sigma_{x_1}, \dots, \Sigma_{x_n}). 
\end{align}  
Next, writing $\Sigma^x: x \mapsto \Sigma(t,x, u)$ [that is, viewing $\Sigma$ as a function of $x$ whilst freezing 
the remaining variables $(t,u)$] we can differentiate \eqref{eq6.44} further to obtain  
\begin{align} \label{eq6.47}
\Delta \Sigma  
=&~\Delta \Sigma^x+e^h \langle \Sigma_{xu} 
+ e^h |\nabla h|^2 \Sigma_u ,\nabla h \rangle \nonumber \\
&+e^h \langle \Sigma_{xu},\nabla h \rangle
+e^{2h} |\nabla h|^2 \Sigma_{uu} +e^h \Sigma_u \Delta h \nonumber\\
& = \Delta \Sigma^x+2e^h \langle \Sigma_{xu} ,\nabla h \rangle 
+ e^h |\nabla h|^2 ( \Sigma_u +e^h \Sigma_{uu})+e^h \Sigma_u \Delta h.
\end{align}   
Next, for $\Delta_f \Sigma$, by using the above calculations and substituting accordingly, we have,  
\begin {align} \label{eq6.48}
\Delta_f \Sigma &= \Delta \Sigma -\langle \nabla f , \nabla \Sigma \rangle = \Delta \Sigma 
-\langle \nabla f , ( \Sigma_x+e^h \Sigma_u\nabla h) \rangle \nonumber\\
& = \Delta \Sigma -\langle \nabla f, \Sigma_x \rangle -e^h \Sigma_u \langle \nabla f , \nabla h \rangle \\
& = \Delta_f \Sigma^x + 2e^h \langle \Sigma_{xu} ,\nabla h \rangle 
+ e^h |\nabla h|^2 ( \Sigma_u +e^h \Sigma_{uu})+e^h \Sigma_u \Delta_f h. \nonumber 
\end{align}
For the sake of future reference we also note that 
\begin {align} \label{eq6.49}
\Delta_f e^{-h} &= \Delta e^{-h} - \langle \nabla f , \nabla e^{-h} \rangle \nonumber \\
&= - {\rm div}(e^{-h}\nabla h)+e^{-h} \langle \nabla f, \nabla h \rangle \nonumber\\
&= -e^{-h} \Delta h + e^{-h} |\nabla h|^2 +e^{-h} \langle \nabla f , \nabla h \rangle 
= -e^{-h} (\Delta_f h - |\nabla h|^2).
\end{align}

Returning to \eqref{1.37} and picking up the estimate from where we left, for the 
last two terms inside the parentheses in the sum on the right, we can write
\begin{align}\label{1.38}
2 (\lambda -1)\langle \nabla h,\nabla (e^{-h} \Sigma) \rangle & + \lambda \Delta_f (e^{-h} \Sigma)  \nonumber\\
=&~2 (\lambda -1) [ - e ^{-h} \Sigma |\nabla h|^2 + e^{-h} \langle \nabla h, \nabla \Sigma \rangle] \nonumber\\
& +\lambda [ e^{-h} \Delta_f  \Sigma + \Sigma \Delta_f e^{-h} -2 e^{-h} \langle \nabla h , \nabla \Sigma \rangle] \nonumber\\
=&-2 (\lambda -1) e ^{-h} \Sigma |\nabla h|^2 + 2 \lambda e^{-h} \langle \nabla h , \nabla \Sigma \rangle \nonumber\\
& -2e^{-h} \langle \nabla h, ( \Sigma_x + e^ h \Sigma_u \nabla h) \rangle 
+ \lambda e^{-h} ( \Delta_f \Sigma^x +2e^h \langle \Sigma_{xu} ,\nabla h \rangle) \nonumber \\
& + \lambda e^{-h} (e^h |\nabla h|^2 ( \Sigma_u +e^h \Sigma_{uu})+e^h \Sigma_u \Delta_f h) \nonumber\\
& + \lambda \Sigma e^{-h} ( -\Delta_f h+ |\nabla h|^2) -2 \lambda e^{-h} \langle \nabla h , \nabla \Sigma \rangle.
\end{align}
As according to \eqref{9.6a} we have 
\begin{align}
\Delta_f h [ \lambda \Sigma_u - \lambda \Sigma e^{-h}] 
&= [-G/(\lambda t_1) - (\lambda -1) |\nabla h|^2/\lambda]
\left[ \lambda(\Sigma_u - \Sigma e^{-h})\right] \nonumber\\
& = -(G/t_1) (\Sigma_u -\Sigma e^{-h}) -(\lambda -1)|\nabla h|^2(\Sigma_u -\Sigma e^{-h}), 
\end{align}  
upon substitution back in \eqref{1.38} this gives 
\begin{align*}
2 (\lambda -1)&\langle \nabla h,\nabla (e^{-h} \Sigma) \rangle +\lambda \Delta_f (e^{-h} \Sigma) \nonumber \\
=&-2 (\lambda -1) e^{-h} \Sigma |\nabla h|^2 
-2 e^{-h}\langle \nabla h, \Sigma_x \rangle  -2 \Sigma_u |\nabla h|^2 \\
& + \lambda e^{-h} \Delta_f \Sigma^x + 2 \lambda \langle \Sigma_{xu} ,\nabla h \rangle 
+ \lambda |\nabla h|^2 \Sigma_u +\lambda |\nabla h|^2 e^h\Sigma_{uu} \\
&+ \lambda \Sigma e^{-h} |\nabla h|^2 - (G/t_1) (\Sigma_u -\Sigma e^{-h})-(\lambda -1)|\nabla h|^2(\Sigma_u - \Sigma e^{-h}) \\
=&~|\nabla h|^2 [-2 (\lambda -1) e^{-h} \Sigma -2 \Sigma_u + \lambda \Sigma_u + 
\lambda e^h \Sigma_{uu} ] \nonumber \\
& + |\nabla h|^2 [ \lambda e^{-h}\Sigma -(\lambda -1) \Sigma_u + (\lambda -1) \Sigma e^{-h}] \nonumber \\
&- [2 \langle \nabla h , (e^{-h}\Sigma_x -\lambda \Sigma_{xu})\rangle
 + (G/t_1) (\Sigma_u -\Sigma e^{-h}) -\lambda e^{-h}\Delta_f \Sigma^x]. 
\end{align*}

Therefore, by taking into account the relevant cancellations, after simplifying terms and using basic inequalities, we can write  
\begin{align} \label{eq6.52}
2 (\lambda -1)&\langle \nabla h, \nabla (e^{-h}\Sigma ) \rangle +\lambda \Delta_f(e^{-h}\Sigma)\nonumber\\
\ge&~|\nabla h|^2(e^{-h}\Sigma - \Sigma_u + \lambda e^h \Sigma_{uu}) - (G/t_1) (\Sigma_u - e^{-h}\Sigma)\nonumber\\
&-2 |\nabla h| |e^{-h} \Sigma_x - \lambda \Sigma_{xu}|+\lambda e^{-h} \Delta_f \Sigma^x.
\end{align}

As a result making use of the relations \eqref{eq6.44}-\eqref{eq6.49} and the inequality \eqref{eq6.52} above and substituting 
all back into \eqref {1.37} and recalling $0 \le \psi \le 1$ we obtain:
\begin{align}\label{1.45}
0 \ge & - \psi G ([c_2 +(m-1) c_1(1+R \sqrt{k_1})+2c_1^2] t_1/R^2 +1+c_1 \underline k_2 t_1) \nonumber \\
& -t_1 \psi^2 G (\Sigma_u - e^{-h}\Sigma) 
+ t_1^2 (\psi ^2/m) [|\nabla h|^2 + e^{-h}\Sigma - \partial_t h]^2  -2c_1 t_1 \psi ^{3/2} |\nabla h| G/R \nonumber\\
& -2 t_1^2 \psi ^2 |\nabla h|^2 [(m-1)k_1 +(\lambda -1) \overline k_2 + k_3 - (e^{-h}\Sigma-\Sigma_u+\lambda e^h \Sigma_{uu})/2] \nonumber\\
& -t_1 ^2 \psi ^2 |\nabla h| (2|e^{-h} \Sigma_x -\lambda \Sigma_{xu}|+ \lambda |\nabla \partial_t f| + 2 \lambda \underline k_2 |\nabla f|) \nonumber\\
&+\lambda t_1^2 \psi^2 (e^{-h}  \Delta_f \Sigma^x-\lambda n [(\underline k_2 + \overline k_2)^2+ 2 k_3]).
\end{align}

In order to obtain the desired bounds out of this it is more efficient to introduce the quantities $y, z$ by setting [as before evaluated at the 
maximum point at $(x_1,t_1)$]
\begin{align} \label{eq6.54}
y= \psi |\nabla h|^2, \qquad z= \psi(\partial_t h -e^{-h}\Sigma). 
\end{align}

Note in particular that $y-\lambda z= \psi G/t_1>0$. Now referring to \eqref{1.45} and recalling the bounds \eqref{k-bounds-two} 
we introduce the constants  
\begin{align} \label{Aeq6.55}
\mathsf{A} = \mathsf{A}^\Sigma =&~2[(m-1)k_1 +(\lambda -1) \overline k_2 + k_3] \nonumber \\
&- \inf_{\Theta_{2R, T}} \{ [(\lambda u^2 \Sigma_{uu} - u\Sigma_u + \Sigma)/u]_-\}, 
\end{align}
and 
\begin{align} \label{Beq6.56}
\mathsf{B} = \mathsf{B}^\Sigma = \lambda \ell_2 + 2 \lambda \underline k_2 \ell_1 
+ \sup_{\Theta_{2R, T}} \{2|(\Sigma_x -\lambda u \Sigma_{xu})/u| \}.  
\end{align}
We remind the reader that here we are making use of the notation introduced earlier on 
\begin{equation}
\Theta_{2R, T} = \{(t,x,u) : (x,t) \in Q_{2R, T}, \, \underline u \le u \le \overline u \} \subset [0, T] \times M \times (0, \infty),
\end{equation} 
where $\overline u$, $\underline u$ denote the maximum and minimum of $u$ on the compact space-time cylinder $Q_{2R, T}$. 
In particular since $u$ is positive we have $[\underline u, \overline u] \subset (0, \infty)$. 
Now substituting the quantities \eqref{eq6.54} back in \eqref{1.45}, recalling again $0 \le \psi \le 1$, 
and utilising \eqref{Aeq6.55}-\eqref{Beq6.56}, basic considerations and bounds lead to   
\begin{align} \label{eq6.51}
0 \ge &  - \psi G ([c_2 +(m-1) c_1(1+R \sqrt{k_1})+2c_1^2] t_1/R^2 +1+c_1 \underline k_2 t_1) \nonumber \\
& + (t_1 ^2/m) [(y-z)^2 - (2mc_1/R) \sqrt y (y-\lambda z)-m\mathsf{A}y -m\mathsf{B} \sqrt y]\nonumber\\
& -t_1 \psi G [ \Sigma_u - e^{-h}\Sigma ]_{+} 
+ \lambda t_1 ^2 \psi ^2 [ e^{-h} \Delta_f \Sigma^x]_- 
- \lambda t_1 ^2 \psi ^2 [\lambda n (\underline k_2 + \overline k_2)^2 + 2 \lambda n k_3]. 
\end{align}

Next an application of Lemma \ref{LiYau-basiclemma} with the choices $a=2mc_1/R$, $b=m\mathsf{A}$ and $c=m\mathsf{B}$ and with 
$y,z$ as in \eqref{eq6.54} and $\lambda>1$ as above gives, for any $\varepsilon \in (0,1)$, 
\begin{align}
(y-z)^2 - & (2m c_1/R) \sqrt y (y-\lambda z) - m\mathsf{A} y -m\mathsf{B} \sqrt y \nonumber\\
\ge&~(y-\lambda z)^2/\lambda^2 - m^2c_1^2 \lambda^2 (y-\lambda z)/[2(\lambda-1)R^2] \nonumber \\
& - m^2 \lambda^2 \mathsf{A}^2 /[4(1-\varepsilon)(\lambda-1)^2]
-(3/4) (m^4 \lambda^2 \mathsf{B}^4/[4 \varepsilon(\lambda-1)^2])^{1/3}. 
\end{align} 
Hence by substituting back in \eqref{eq6.51} it follows that   
\begin{align} \label{eq6.55}
0 \ge &~- \psi G ([c_2 +(m-1) c_1(1+R \sqrt{k_1})+2c_1^2] t_1/R^2 +1+c_1 \underline k_2 t_1)  \nonumber\\
&+ (t_1^2/m) [(\psi G)^2/(t_1^2 \lambda ^2) 
- m^2 c_1 ^2 \lambda^2 (\psi G) / (2(\lambda -1)R^2 t_1)] \nonumber \\
& -(mt_1^2\lambda^2 \mathsf{A}^2) /[4(1-\varepsilon)(\lambda-1)^2] -t_1 \psi G [\Sigma_u - e^{-h}\Sigma]_{+} \nonumber \\
& -[(3t_1^2)/(4m)] (m^4 \lambda^2 \mathsf{B}^4/[4 \varepsilon(\lambda-1)^2])^{1/3} \nonumber \\
&+ \lambda t_1 ^2 \psi ^2 [ e^{-h} \Delta_f \Sigma^x ]_- 
- \lambda t_1 ^2 \psi ^2 [\lambda n (\underline k_2 + \overline k_2)^2 + 2 \lambda n k_3].
\end{align}
Upon setting 
\begin{align} \label{Deq6.60}
\mathsf{D} =&~[c_2 +(m-1) c_1(1+R \sqrt{k_1})+2c_1^2]t_1/R^2 +1 \nonumber \\
&+ c_1 \underline k_2 t_1 + m t_1 c_1 ^2 \lambda ^2/[2(\lambda -1)R^2] + t_1 \gamma^\Sigma_1, 
\end{align}
and 
\begin{align} \label{Eeq6.61}
\mathsf{E}=&~m \lambda^2 \mathsf{A}^2/[4(1-\varepsilon)(\lambda-1)^2] \nonumber \\
&+ (3/4) [ m \lambda^2 \mathsf{B}^4/(4 \varepsilon(\lambda-1)^2) ]^{1/3} \nonumber \\
&+ \lambda ^2 n (\underline k_2 + \overline k_2)^2 + 2 \lambda^2 n k_3 + \lambda \gamma^\Sigma_2, 
\end{align}
where 
\begin{equation} \label{gamma-one-eq6.64}
\gamma^\Sigma_1 = \sup_{\Theta_{2R, T}} \{[(u\Sigma_u - \Sigma)/u]_{+}\}, \qquad 
\gamma^\Sigma_2 = - \inf_{\Theta_{2R, T}} \{[\Delta_f \Sigma^x/u]_-\}, 
\end{equation}
we can rewrite \eqref{eq6.55} after rearranging terms as 
\begin{align} \label{eq6.62}
0 \ge (\psi G)^2/(m\lambda^2) - (\psi G) \mathsf{D} - t_1^2 \mathsf{E}.
\end{align}
As a result basic considerations on the inequality \eqref{eq6.62} lead to the conclusion 
\begin {align} \label{eq6.63}
\psi G &\le (m \lambda^2/2) \left(  \mathsf{D} + \sqrt{\mathsf{D}^2+ 4 t^2_1 \mathsf{E}/(m\lambda^2)} \right) \nonumber\\
&\le (m \lambda^2/2) \left( 2\mathsf{D}+\sqrt{(4t_1^2 \mathsf{E})/(m \lambda^2)} \right) \nonumber \\
&= m \lambda^2\mathsf{D} +t_1 \lambda \sqrt{m\mathsf{E}}.
\end{align}
Since $\psi \equiv 1$ for $d(x,x_0, \tau) \le R$ and $(x_1, t_1)$ is the point where $\psi G$ attains its maximum 
on $\{d(x,x_0, t) \le 2R, 0 \le t \le \tau\}$ we have
\begin{align}
G(x, \tau) = [\psi G] (x, \tau) \le [\psi G] (x_1,t_1) \le m \lambda^2\mathsf{D} +t_1 \lambda \sqrt{m\mathsf{E}}.
\end{align}
Therefore recalling \eqref{2.3}, substituting for $\mathsf{D}$ and $\mathsf{E}$ from \eqref{Deq6.60} and \eqref{Eeq6.61} above and 
making noting $t_1 \le \tau$, we can write after dividing both sides $\lambda \tau$, 
\begin{align}
\lambda^{-1} |\nabla h|^2 - \partial_t h + e^{-h} \Sigma
\le&~(m \lambda/\tau) \mathsf{D} + \sqrt{m \mathsf{E}}  \nonumber \\
\le&~(m \lambda) [c_2 +(m-1) c_1(1+R \sqrt{k_1})+2c_1^2] /R^2 \nonumber \\
& + (m \lambda/\tau) 
+ (m \lambda) (\gamma^\Sigma_1 + c_1 \underline k_2 + m c_1 ^2 \lambda ^2/[2(\lambda -1)R^2])  \nonumber \\
& + \sqrt m \{ m \lambda^2 \mathsf{A}^2/[4(1-\varepsilon)(\lambda-1)^2] \nonumber \\
& + (3/4) [ m \lambda^2 \mathsf{B}^4/(4 \varepsilon(\lambda-1)^2) ]^{1/3} \nonumber \\
& + \lambda ^2 n (\underline k_2 + \overline k_2)^2 + 2 \lambda^2 n k_3 + \lambda \gamma^\Sigma_2 \}^{1/2}. 
\end{align}
Finally using the arbitrariness of $0< \tau \le T$ it follows after reverting back to $u$ upon noting the relation $h=\log u$ 
and rearranging terms that 
\begin{align}
\frac{|\nabla u|^2}{\lambda u^2} - \frac{\partial_t u}{u} + \frac{\Sigma}{u} 
\le&~(m \lambda) [1/t + \gamma^\Sigma_1 + c_1 \underline k_2] \nonumber \\
&+ (m \lambda) [m c_1 ^2 \lambda ^2/[2(\lambda -1)]+c_2 +(m-1) c_1(1+R \sqrt{k_1})+2c_1^2]/R^2 \nonumber \\
&+ \sqrt m \{ m \lambda^2 \mathsf{A}^2/[4(1-\varepsilon)(\lambda-1)^2]
+ (3/4) [ m \lambda^2 \mathsf{B}^4/(4 \varepsilon(\lambda-1)^2) ]^{1/3}  \nonumber \\
&+ \lambda ^2 n (\underline k_2 + \overline k_2)^2 + 2 \lambda^2 n k_3 
+ \lambda \gamma^\Sigma_2 \}^{1/2}. 
\end{align}
which is the desired estimate \eqref{1.26} with $q \equiv 0$. Now to establish the estimate in its full strength (for general $q$) 
it suffices to replace $\Sigma$ with $\Sigma+qu$ and use the conclusion from the first part. Then making note of 
[{\it cf.} \eqref{eq2.17} and \eqref{gamma-one-eq6.64}]  
\begin{align}
\gamma_1^{\Sigma + qu} = \gamma_1^\Sigma, \qquad 
\gamma_2^{\Sigma + qu} \le \gamma_2^\Sigma + \gamma_2^{qu}, 
\end{align}
along with $\mathsf{A}^{\Sigma+qu}=\mathsf{A}^\Sigma$ [{\it cf.} \eqref{eq2.13} and \eqref{Aeq6.55}] and 
$\mathsf{B}$ in \eqref{Beq6.56} changing to \eqref{eq2.15} gives the estimate \eqref{1.26}. The proof is thus 
complete. \hfill $\square$

\section{Proof of the parabolic Harnack inequality in Theorem \ref{thm38}} \label{sec7}

With the aid of the estimates established in Theorem \ref{thm28} we can now prove the desired parabolic 
Harnack inequality in Theorem \ref{thm38}. Towards this end it suffices to integrate the former estimate 
along suitable space-times curves in $Q_{R,T} \subset M \times [0, T]$. Here we shall prove only the 
local Harnack inequality. The global inequality is similar (see the comments at the end).
Towards this end let us first move on to rewriting the Li-Yau Harnack inequality \eqref{1.26} as follows:
\begin{align} \label{Heq8.1}
\frac{|\nabla u|^2}{\lambda u^2} - \frac{\partial_t u}{u} + q+ \frac{\Sigma(t,x,u)}{u} 
\le&~\frac{m \lambda}{t} + m \lambda (\gamma^\Sigma_1 + c_1 \underline k_2) \nonumber \\
&+ \frac{m \lambda}{R^2} \left[ \frac{m c_1 ^2 \lambda ^2}{2(\lambda -1)}+c_2 +(m-1) c_1(1+R \sqrt{k_1})+2c_1^2 \right]  \nonumber \\
&+ \sqrt m \bigg\{ \frac{m \lambda^2 \mathsf{A}^2}{4(1-\varepsilon)(\lambda-1)^2} 
+ \frac{3}{4} \left[ \frac{m \lambda^2 \mathsf{B}^4}{4 \varepsilon(\lambda-1)^2} \right]^{1/3} \nonumber \\
&+ \lambda ^2 n (\underline k_2 + \overline k_2)^2 + 2 \lambda^2 n k_3 
+ \lambda (\gamma^\Sigma_2 + \gamma_2^{qu}) \bigg\}^{1/2}.   
\end{align}
Put $\overrightarrow{k} = (k_1, \underline k_2, \overline k_2, k_3)$, 
$\overrightarrow{\gamma}=(\gamma^\Sigma_1, \gamma^\Sigma_2, \gamma_2^{qu}, \gamma^\Sigma_3)$
and let 
$\mathsf{S}=\mathsf{S}(m, \varepsilon, \lambda, \mathsf{A}, \mathsf{B}, R, T, \overrightarrow{k}, \overrightarrow{\gamma}, \underline q)$ 
be defined by 
\begin{align} \label{S-eq7.1}
\mathsf{S} =&~
- \frac{m \lambda}{R^2} \left[ \frac{m c_1 ^2 \lambda ^2}{2(\lambda -1)}+c_2 +(m-1) c_1(1+R \sqrt{k_1})+2c_1^2 \right] \nonumber \\
&~-\sqrt m \bigg\{ \frac{m \lambda^2 \mathsf{A}^2}{4(1-\varepsilon)(\lambda-1)^2} 
+ \frac{3}{4} \left[ \frac{m \lambda^2 \mathsf{B}^4}{4 \varepsilon(\lambda-1)^2} \right]^{1/3}  \nonumber \\
&~+\lambda ^2 n (\underline k_2 + \overline k_2)^2 + 2 \lambda^2 n k_3 
+ \lambda (\gamma^\Sigma_2 + \gamma^{qu}_2) \bigg\}^{1/2} \nonumber \\
&~-m \lambda (\gamma^\Sigma_1 + c_1 \underline k_2) + \underline q + \gamma^\Sigma_3,  
\end{align}
where 
\begin{equation}
\gamma^\Sigma_3 = \inf_{\Theta_{2R,T}} \left\{ \frac{\Sigma(t,x,u)}{u} \right\}, \qquad \underline q = \inf_{Q_{2R,T}} q = \gamma_3^{qu}.
\end{equation}
It follows from \eqref{Heq8.1} that $\partial_t u/u \ge |\nabla u|^2/(\lambda u^2) - m \lambda/t + \mathsf{S}$. 
Suppose $\gamma \in \mathscr{C}^1( [t_1,t_2]; M)$ is an arbitrary curve lying entirely in $Q_{R, T}$ with $\gamma(t_1) = x_1$ 
and $\gamma(t_2) = x_2$. Using the above it is seen that 
\begin{align}
d/dt [\log u(\gamma(t),t)] &= \langle \nabla u/u, \dot \gamma (t) \rangle + \partial_t u/u \nonumber \\
&\ge  \langle \nabla u/u, \dot \gamma (t) \rangle + |\nabla u|^2/(\lambda u^2) - ( m \lambda)/t + \mathsf{S}\nonumber\\
&= \lambda^{-1} |\nabla u/u + \lambda \dot \gamma (t)/2|^2 - \lambda |\dot \gamma (t)|^2/4 
-(m \lambda)/t + \mathsf{S} \nonumber\\ 
&\ge -\lambda |\dot \gamma (t)|^2/4 - (m \lambda)/t + \mathsf{S}, 
\end{align}
where the inner products are with respect to $g(t)$. Therefore integrating the above inequality gives 
\begin{align}
\log \frac{u(x_2,t_2)}{u(x_1,t_1)} &= \log u(\gamma(t),t)\bigg|_{t_1}^{t_2}= \int_{t_1}^{t_2} \frac{d}{dt} \log u(\gamma(t),t) \, dt \nonumber \\
&\ge \int_{t_1}^{t_2}  -\frac{\lambda}{4} |\dot \gamma (t)|^2 \,dt 
- \int_{t_1}^{t_2} \frac{m \lambda}{t} dt+\int_{t_1}^{t_2} \mathsf{S} \, dt \nonumber \\
&= -m \lambda \log (t_2/t_1) - (\lambda/4) \int_{t_1}^{t_2} |\dot \gamma (t)|^2 \,dt + (t_2-t_1) \mathsf{S}.
\end{align}
Hence upon exponentiating we have
\begin{align}
\frac{u(x_2,t_2)}{u(x_1,t_1)} &\ge \left(\frac{t_2}{t_1} \right)^{-m \lambda} 
\exp \left[ - \int_{t_1}^{t_2} \frac{\lambda}{4}|\dot \gamma(t)|^2 \, dt \right]  \exp [(t_2-t_1) \mathsf{S}]  ,  
\end{align}
or upon rearranging terms and rescaling the integral:  
\begin{align}
u(x_2,t_2)\geq u(x_1,t_1)\left(\frac{t_2}{t_1}\right)^{-m \lambda}e^{- \lambda L(x_1,x_2, t_2-t_1)} e^{(t_2-t_1) \mathsf{S}} 
\end{align}
where 
\begin{align}
L(x_1,x_2, t_2-t_1)=\inf_\gamma \left[ \frac{1}{4(t_2-t_1)} \int_{0}^{1} |\dot \gamma(t)|^2\,dt \right].
\end{align}

This gives the parabolic Harnack inequality in its local form. Now if the bounds are global by arguing exactly 
as above using the global estimate in Theorem \ref{thm28-global} we obtain the global counterpart of the 
inequality. This therefore completes the proof. \hfill $\square$

\section{Proof of the Liouville results in Theorem \ref{thm46} and Theorem \ref{thm48}} \label{sec8}

\qquad \\
{\bf Proof of Theorem \ref{thm46}.} 
Under the stated assumptions and the fact that $u$ and the metric-potential pair $(g,f)$ are time independent 
(i.e., $\partial_t u \equiv 0$, $\partial_t g \equiv 0$, $\partial_t f \equiv 0$), it follows from the global estimate 
\eqref{eq1.12-global} in Theorem \ref{thm18-global} with $q \equiv 0$ and $k=0$ that 
\begin{align} \label{eqL8.1}
\sup_M \left( \frac{|\nabla u|}{\sqrt u} \right) 
\le & C \sup_{M} \left( \mathsf{T}_\Sigma^{1/2} (u)+ \mathsf{S}_\Sigma^{1/3} (u) \right) 
\Big( \sup_{M} \sqrt u \Big), 
\end{align}
where $\mathsf{T}_\Sigma (u)$ and $\mathsf{S}_\Sigma(u)$ are the expressions given in \eqref{eq2.8}. 
Now a close inspection of the latter expressions upon recalling the condition $\Sigma=\Sigma(u)$ gives $\mathsf{S}_\Sigma(u)=0$ 
while by virtue of the inequality $\Sigma(u)-2u \Sigma_u(u) \ge 0$ we have
\begin{equation} 
\mathsf{T}_\Sigma (u) = \left[ \frac{2u \Sigma_u(u)-\Sigma(u)}{u} \right]_+ =0.  
\end{equation}
As a result it follows from \eqref{eqL8.1} and the global bound $\sup_M \sqrt u = \sqrt{\sup_M u} < \infty$ that $|\nabla u| \equiv 0$ and so the conclusion follows. 
\hfill $\square$

\qquad \\
{\bf Proof of Theorem \ref{thm48}.} 
The proof of \eqref{eqL2.26} follows directly from the global estimate \eqref{1.26-global} in Theorem \ref{thm28-global} 
upon noting that $u$ and the metric-potential pair $(g,f)$ are time independent 
(i.e., $\partial_t u \equiv 0$, $\partial_t g \equiv 0$, $\partial_t f \equiv 0$). Now combining the latter with the assumptions 
on the nonlinearity $\Sigma$ and the solution $u$ as formulated in the theorem it follows that   
\begin{align} 
\frac{|\nabla u|^2}{\lambda u^2} + \frac{\Sigma(u)}{u} 
\le&~ m \lambda 
\sup_{\Theta} \left\{ \left[\frac{u \Sigma_u (u) - \Sigma(u)}{u}\right]_{+} \right\} \nonumber \\
&+ \frac{m \lambda/\sqrt{1-\varepsilon}}{2(\lambda-1)} 
\sup_{\Theta} \left\{ \left[ \frac{-\Sigma(u) + u\Sigma_u(u) - \lambda u^2 \Sigma_{uu}(u)}{u} \right]_+ \right\} =0.  
\end{align}
Since $\Sigma(u) \ge 0$ we thus infer that $|\nabla u|^2/(\lambda u^2) + \Sigma(u)/u \equiv 0$ 
and therefore $|\nabla u| \equiv 0$ on $M$. The conclusion on $u$ being a constant is now immediate. 
\hfill $\square$

\section{Hamilton type estimates and universal global bounds on closed manifolds} \label{sec9}

\begin{lemma}\label{lemma31}
Let u be a bounded positive solution to $(\partial_t - q(x,t) -\Delta_f) u = \Sigma(t,x,u)$ with $0<u\leq D$ 
and suppose the metric-potential pair $(g,f)$ is of class $\mathscr{C}^2$ and evolves under the 
$\mathsf{k}$-super Perelman-Ricci flow \eqref{eq11c}. Let
\begin{align} \label{eq12.1}
\mathscr{F}_\gamma[u]  = \gamma(t) |\nabla u|^2/u -u \log(D/u), 
\end{align}
where $\gamma$ is a smooth, non-negative but otherwise arbitrary function. Then we have
\begin{align} \label{eq12.11}
\square_q \mathscr{F}_\gamma [u] 
\le&~2 \gamma \langle \nabla u ,\nabla q \rangle + (\gamma^{'} + 2 \mathsf{k} \gamma -1) |\nabla u|^2/u 
+ 2(\gamma/u) \langle \nabla u , \nabla \Sigma \rangle \nonumber\\
&+[1-\log(D/u)] \Sigma(t,x,u) - \gamma (|\nabla u|^2/u^2) \Sigma(t,x,u) + qu, 
\end{align}
where as before $\square_q = (\partial_t - q(x,t) -\Delta_f)$. 
\end {lemma}

\begin{proof}
Starting from \eqref{eq12.1} and working our way forward, it is firstly seen that, 
\begin{align}
\left[\begin{array}{c} \nabla \\ \partial_t \end{array}\right] 
\frac{|\nabla u|^2}{u} = 
\left[\begin{array}{c} (1/u) \nabla |\nabla u|^2 - (|\nabla u|^2/u^2) \nabla u \\ 
(2/u) \langle \nabla u, \nabla \partial_t u \rangle - (|\nabla u|^2/u^2) \partial_t u - (1/u) [\partial_t g](\nabla u, \nabla u)
\end{array}\right].   
\end{align}
A straightforward calculation upon further differentiation and forming the expression $\Delta_f (|\nabla u|^2/u)$ then gives 
\begin{align}
\Delta_f \frac{|\nabla u|^2}{u} 
&= \Delta \frac{|\nabla u|^2}{u} - \left\langle \nabla f, \nabla \frac{|\nabla u|^2}{u} \right\rangle \nonumber \\ 
&= \frac{1}{u} \Delta_f |\nabla u|^2 - 2\frac{\langle \nabla |\nabla u|^2 , \nabla u \rangle }{u^2}
 -\frac{|\nabla u|^2}{u^2} \Delta_f u + 2\frac{|\nabla u|^4}{u^3}.
\end{align}
Similar calculations for the second term $u \log (D/u)$ in \eqref{eq12.1} lead to   
\begin{align}
\left[\begin{array}{c} \nabla \\ \partial_t \end{array}\right] 
[u \log (D/u)]
= [\log (D/u) -1]
\left[\begin{array}{c} \nabla \\ \partial_t \end{array}\right] u, 
\end{align}
and so after a further differentiation result in  
\begin{align}
\Delta_f [u \log (D/u)] = [\log (D/u) -1] \Delta_f u - |\nabla u|^2/u.  
\end{align}

An application of the weighted $q$-heat operator $\square_q= (\partial_t - q -\Delta_f)$ to $\mathscr{F}_\gamma [u]$ 
in \eqref{eq12.1} and putting together the above fragments give  
\begin{align}
\square_q \mathscr{F}_\gamma [u]  
=&~\gamma ^{'} [|\nabla u|^2/u] + \gamma (\partial_t - q-\Delta_f) [|\nabla u|^2/u]  
- (\partial_t - q-\Delta_f) [u \log (D/u)] \nonumber\\
=&~ \gamma ^{'} \frac{|\nabla u|^2}{u} +\frac{\gamma}{u} \left[ - \partial_t g (\nabla u,\nabla u) 
+ 2 \langle \nabla u , \nabla \partial_t u \rangle - \frac{|\nabla u|^2 }{u} \partial_t u \right] \nonumber\\
& - \gamma q \frac{|\nabla u|^2}{u} - \gamma \left[ \frac{\Delta_f |\nabla u|^2 }{u} 
- \frac{2 \langle \nabla |\nabla u|^2 , \nabla u \rangle }{u^2} 
- \Delta_f u \frac{|\nabla u|^2}{u^2} + 2 \frac{|\nabla u |^4}{u^3} \right] \nonumber\\
& - [\log (D/u) -1] \partial_t u + q u \log (D/u) + [\log (D/u) -1] \Delta_f u -\frac{|\nabla u|^2}{u} \nonumber\\
& - qu +qu - \gamma q \frac{|\nabla u|^2}{u} + \gamma q \frac{|\nabla u|^2}{u}, 
\end{align}
or after some calculation and rearrangement of terms 
\begin{align} \label{4.10}
\square_q \mathscr{F}_\gamma [u]  
=&~(\gamma ^{'} -1) \frac{|\nabla u|^2}{u} + \frac{\gamma}{u} [ - \partial_t g (\nabla u,\nabla u) 
+ 2 \langle \nabla u , \nabla \partial_t u \rangle - \Delta_f |\nabla u|^2  \nonumber\\
&+ \frac{2}{u} \langle \nabla |\nabla u|^2 , \nabla u \rangle -2 \frac{|\nabla u|^4}{u^2}] 
- [\log \frac{D}{u} -1](\partial_t - q -\Delta_f) u\nonumber\\
&-\gamma \frac{|\nabla u|^2}{u^2} (\partial_t - q -\Delta_f) u 
+ qu -2 \gamma q \frac{|\nabla u|^2}{u}.
\end{align}
Now by recalling the equation $\square_q u = (\partial_t - q -\Delta_f) u = \Sigma (t,x,u)$ satsfied by $u$ and hence that 
$2 \langle \nabla u , \nabla \partial_t u \rangle 
= 2 \langle \nabla u , \nabla \Delta_f u \rangle + 2 \langle \nabla u , \nabla (qu) \rangle + 2 \langle \nabla u, \nabla \Sigma \rangle$ 
and upon invoking the Bochner-Weitzenb\"ock formula 
$\Delta_f |\nabla u|^2 = 2 |\nabla ^2 u|^2 + 2\langle \nabla u , \nabla \Delta_f u \rangle + 2 {\mathscr Ric}_f (\nabla u, \nabla u)$ 
we can rewrite \eqref{4.10} after substitution from the above as 
\begin{align}
\square_q \mathscr{F}_\gamma [u] 
=&~(\gamma ^{'} -1) \frac{|\nabla u|^2}{u} - \frac{\gamma}{u} [ \partial_t g (\nabla u,\nabla u) +2 |\nabla ^2 u|^2+ 
2 {\mathscr Ric}_f (\nabla u, \nabla u)\nonumber\\
& - \frac{2}{u} \langle \nabla |\nabla u|^2 , \nabla u \rangle] - 2 \frac{\gamma }{u} \frac{|\nabla u|^4}{u^2} 
+ 2 \frac{\gamma}{u} \langle \nabla u, \nabla (qu) \rangle 
+ 2 \frac{\gamma}{u} \langle \nabla u, \nabla \Sigma \rangle \nonumber\\
& - [\log (D/u) -1] \Sigma(t,x,u) -\gamma \frac{|\nabla u|^2}{u^2} [ \Sigma ( t,x,u) + 2qu] + qu.
\end{align}
Using basic tensor algebra and making note of the the non-negativity of the expression   
\begin{align}
|\nabla ^2 u|^2 - \frac{\langle \nabla |\nabla u|^2 , \nabla u \rangle }{u} + \frac{|\nabla u|^4}{u^2} 
= \left| \nabla ^2 u - \frac{\nabla u \otimes \nabla u }{u} \right| ^2 \ge 0, 
\end{align}
followed by an application of the Perelman-Ricci flow inequality $\partial_t g + 2 {\mathscr Ric}_f (g) \ge -2 \mathsf{k} g$ 
as satisfied by $g$, we arrive at  
\begin{align}
\square_q \mathscr{F}_\gamma [u] 
\le&~2 \gamma \langle \nabla u ,\nabla q \rangle + (\gamma^{'} + 2 \mathsf{k} \gamma -1)\frac{|\nabla u|^2}{u} 
+ 2\frac{\gamma}{u} \langle \nabla u , \nabla \Sigma \rangle \nonumber \\ 
&+ [1-\log(D/u)]\Sigma(t,x,u) - \gamma \frac{|\nabla u|^2}{u^2} \Sigma(t,x,u) + qu, 
\end{align}
which is the required conclusion.
\end{proof}

\qquad \\
{\bf Proof of Theorem \ref{thm25}.} 
The function $\gamma(t)=t/(1+2\mathsf{k}t)$ is non-negative and satisfies $\gamma'+2\mathsf{k} \gamma-1 \le 0$. 
Applying Lemma \ref{lemma31} with $q=0$ and $eD$ in place of $D$ (note that  
$u \le D \implies u \le eD$) we have from \eqref{eq12.11}
\begin{align}
\square_q \mathscr{F}_\gamma [u] 
\le&~2 \gamma \langle \nabla u ,\nabla q \rangle + (\gamma^{'} + 2\mathsf{k} \gamma -1)\frac{|\nabla u|^2}{u} 
+ 2\frac{\gamma}{u} \langle \nabla u , \nabla \Sigma(t,x,u) \rangle \nonumber \\ 
&+ [1-\log(e D/u)]\Sigma(t,x,u) - \gamma \frac{|\nabla u|^2}{u^2} \Sigma(t,x,u) + qu \nonumber \\
\le&~(\gamma^{'} + 2\mathsf{k} \gamma -1)\frac{|\nabla u|^2}{u} 
+ 2\frac{\gamma}{u} \langle \nabla u, \Sigma_x(t,x,u) \rangle \nonumber \\ 
&~+ 2\frac{\gamma}{u} |\nabla u|^2 \Sigma_u(t,x,u) 
- [\log(D/u)] \Sigma(t,x,u) - \gamma \frac{|\nabla u|^2}{u^2} \Sigma(t,x,u) \nonumber \\
\le &~2\frac{\gamma |\nabla u|^2}{u} \Sigma' (u) 
- \frac{\gamma |\nabla u|^2}{u^2} \Sigma(u) - \log (D/u) \Sigma(u) \le 0.\nonumber 
\end{align}

Next by an easy inspection $\mathscr{F}_\gamma[u](x,0) \le 0$ for all $x \in M$. Furthermore as seen 
above $\square \mathscr{F}_\gamma[u] = (\partial_t - \Delta_f) \mathscr{F}_\gamma[u] \le 0$. Thus, 
since $M$ is closed, an application of maximum principle gives $\mathscr{F}_\gamma[u](x,t)  \le 0$ 
for all $(x, t)$ in $M \times [0, T]$ from which the desired  estimate \eqref{eq15} follows. 
Next to prove \eqref{eq16} set $\mathsf{U}(x,t) = \log [e D/u(x,t)]$. Then by (\ref{eq15}),  
\begin{equation}
\left| \nabla \sqrt{\mathsf{U}} \right| = \left| \frac{\nabla u/u}{\sqrt{4\mathsf{U}}} \right| 
\le \frac{\sqrt{1+2\mathsf{k}t}}{\sqrt{4t}}.
\end{equation}
Integrating the above along a minimising geodesic joining a pair of points $x_1,  x_2$ in $M$ then gives 
\begin{equation}
\sqrt{\log (e D/u(y,t))} - \sqrt{\log (e D/u(x,t))} \le d(x,y; t) \frac{\sqrt{1+2\mathsf{k}t}}{\sqrt{4t}}.
\end{equation} 
For any $s>0$ thus 
$\log (e D/u(y,t)) \leq (1+s) [ \log (e D/u(x,t)) + d^2(x,y;t) (1+2\mathsf{k}t)/(4st) ]$. 
Exponentiating and rearranging yields the desired inequality \eqref{eq16}. 
\hfill $\square$


\qquad \\
{\bf Acknowledgement.} The authors gratefully acknowledge support from EPSRC.

\end{document}